%% file: main.tex
\documentclass[a4paper,fleqn]{cas-sc}

\usepackage[authoryear,longnamesfirst]{natbib}

\usepackage[utf8]{inputenc}
\usepackage[british]{babel}
\usepackage{amsmath,amsfonts,amssymb,bm,amsthm}
\usepackage{numprint}

\usepackage{algorithm}
\usepackage{algpseudocode}
\usepackage{float}
\usepackage{tikz}
\usepackage{forest}
\usetikzlibrary{arrows}
\usetikzlibrary {arrows.meta} 
\usetikzlibrary{shadows}
\usetikzlibrary{calc}

\usepackage{color}
\definecolor{imperialblue}{RGB}{0 62 116}
\definecolor{lightblue}{RGB}{212 239 252}
\definecolor{lightgrey}{RGB}{235 238 238}
\definecolor{navy}{RGB}{0 33 71}
\definecolor{blue}{RGB}{0 110 175}
 \definecolor{processblue}{RGB}{0 133 202}
\definecolor{poolblue}{RGB}{12 161 205}
\definecolor{darkteal}{RGB}{15 130 145}
\definecolor{seaglass}{RGB}{55 159 159}
\definecolor{darkgreen}{RGB}{2 137 59}
\definecolor{kermitgreen}{RGB}{102 164 10}
\definecolor{lime}{RGB}{169 214 0}
\definecolor{cyan}{RGB}{0 255 255}
\definecolor{tangerine}{RGB}{236 115 0}
\definecolor{violet}{RGB}{150 0 120}
\definecolor{red}{RGB}{221 37 1}

\colorlet{mygreen}{green!25!black}
\colorlet{col1in}{red!30}
\colorlet{col1out}{red!40}
\colorlet{col2in}{mygreen!40}
\colorlet{col2out}{mygreen!50}
\colorlet{col3in}{blue!30}
\colorlet{col3out}{blue!40}
\colorlet{col4in}{mygreen!20}
\colorlet{col4out}{mygreen!30}
\colorlet{col5in}{blue!10}
\colorlet{col5out}{blue!20}
\colorlet{col6in}{blue!20}
\colorlet{col6out}{blue!30}
\colorlet{col7out}{orange}
\colorlet{col7in}{orange!50}
\colorlet{col8out}{orange!40}
\colorlet{col8in}{orange!20}
\colorlet{linecol}{blue!60}

\pgfkeys{/forest,
        rect/.append style={rectangle, rounded corners=2pt, inner color=col6in, outer color=col6out},
        rrect/.append style={rectangle, rounded corners=4pt, inner color=col6in, outer color=col6out},
        ellip/.append style={ellipse, inner color=col5in, outer color=col5out},
        orect/.append style={rect, font=\sffamily\bfseries\LARGE, text width=325pt, text centered, minimum height=10pt, outer color=col7out, inner color=col7in},
        oellip/.append style={ellip, inner color=col8in, outer color=col8out, font=\sffamily\bfseries\large, text centered},
}

\usepackage{xcolor}
\usepackage{listings} % Display code / shell commands
\usepackage{lstautogobble}
\lstdefinestyle{C++Style}{
  language=C++,
  basicstyle=\ttfamily\footnotesize,
  numbers=left,
  numberstyle=\tiny,
  numbersep=3pt,
  frame=,
  columns=fullflexible,
  backgroundcolor=\color{black!05},
  linewidth=0.95\linewidth,
  xleftmargin=0.1\linewidth,
  showspaces=false,
  showstringspaces=false,
  breaklines=true,
  keywordstyle=\color{blue}\ttfamily,
  stringstyle=\color{red}\ttfamily,
  commentstyle=\color{violet}\ttfamily,
  morecomment=[l][\color{teal}]{\#}
}
\lstdefinelanguage{XML}
{
  basicstyle=\ttfamily\footnotesize,
  morestring=[b]",
  moredelim=[s][\bfseries\color{maroon}]{<}{\ },
  moredelim=[s][\bfseries\color{maroon}]{</}{>},
  moredelim=[l][\bfseries\color{maroon}]{/>},
  moredelim=[l][\bfseries\color{maroon}]{>},
  morecomment=[s]{<?}{?>},
  morecomment=[s]{<!--}{-->},
  commentstyle=\color{gray},
  stringstyle=\color{orange},
  identifierstyle=\color{darkblue},
  showstringspaces=false
}
\lstdefinestyle{XMLStyle}{
  language=XML,
  basicstyle=\ttfamily\footnotesize,
  numbers=left,
  numberstyle=\tiny,
  numbersep=3pt,
  frame=,
  columns=fullflexible,
  backgroundcolor=\color{black!05},
  linewidth=\linewidth,
  xleftmargin=0.05\linewidth,
  keepspaces=true
}

\newcommand{\dd}{\,\textrm d}   % infintesimal

\algrenewcommand\algorithmicrequire{\textbf{Input:}}
\algrenewcommand\algorithmicensure{\textbf{Output:}}

\usepackage{url}

\def\tsc#1{\csdef{#1}{\textsc{\lowercase{#1}}\xspace}}
\tsc{WGM}
\tsc{QE}
\tsc{EP}
\tsc{PMS}
\tsc{BEC}
\tsc{DE}
\begin{document}
\let\WriteBookmarks\relax
\def\floatpagepagefraction{1}
\def\textpagefraction{.001}
\shorttitle{Architecture-aware $h$-to-$p$ optimisation}
\shortauthors{J.Y. Xing et~al.}

\title[mode = title]{Architecture-aware $h$-to-$p$ optimisation: spectral/$hp$ element operators for mixed-element meshes}

\author[1]{Jacques Y. Xing}
\cormark[1]
\ead{j.xing@imperial.ac.uk}

\author[2]{Boyang Xia}
\cormark[1]
\ead{boyang.xia@kcl.ac.uk}

\author[1]{Diego Renner}

\author[1]{Chris D. Cantwell}

\author[2] {David Moxey}

\author[3] {Robert M. Kirby}

\author[1] {Spencer J. Sherwin}

\affiliation[1]{organization={Imperial College London},
                addressline={Exhibition Rd, South Kensington}, 
                city={London},
                postcode={SW7 2AZ}, 
                state={England},
                country={United Kingdom}}
                
\affiliation[2]{organization={King's College London},
                addressline={Strand}, 
                postcode={WC2R 2LS}, 
                city={London},
                state={England},
                country={United Kingdom}}
                
\affiliation[3]{organization={University of Utah},
                addressline={201 Presidents' Cir}, 
                postcode={UT 84112}, 
                city={Salt Lake City},
                state={Utah},
                country={United States}}

\cortext[cor1]{Lead Co-authors}

\begin{abstract}
We extend earlier international efforts to optimise hexahedral-based spectral element methods on GPUs and vectorised CPUs to mixed element meshes additionally involving prismatic, pyramidic, and tetrahedral shapes using tensorial expansions. We demonstrate that common finite element operators (such as the mass and Helmholtz matrices) benefit from alternative implementation strategies depending on the element shape, choice of polynomial order, and system architecture in order to achieve optimal performance. In addition, we introduce a new approach/interpretation to efficiently evaluate more complex operations involving inner products with the derivative of the expansions as part of the integrand such as the stiffness matrix. This approach seeks to maximise operations using the collocation properties of the nodal tensorial expansion associated with classical quadrature rules. Our GPU performance tests demonstrate that the throughput of the Helmholtz operator on tetrahedral elements is at most $2.5$ times slower than on hexahedral elements, despite tetrahedra having a factor of six greater floating-point operations. 
\end{abstract}


\maketitle

\input{Introduction}
\input{Theory}
\input{Implementation}
\input{Results}
\input{Summary}

\input{Appendix}

\bibliographystyle{unsrt}
\bibliography{refs}

\vskip3pc

\end{document}

%% file: Introduction.tex
\section{Introduction}
High-order finite elements are increasingly being used to solve partial differential equations (PDEs) to model sophisticated physics in geometrically complex domains. There is a growing interest in their industrial application, particularly in the area of fluid dynamics where their beneficial numerical properties, such as low numerical diffusion and dispersion, allow transient dynamics to be accurately captured over long integration times. The high-order nature of the finite element expansion bases also readily aligns with modern architectures due to their higher arithmetic intensity, by allowing more floating-point operations (FLOPS) to be performed on each byte of data transferred onto the processor. This has also made these methods of particular interest for use on GPUs.

Nektar++ \cite{nektar2015} is traditionally a CPU-based framework for solving PDEs using spectral/$hp$ element methods on mixed-element meshes. The support for simplex elements enables effective use of the tools on challenging industrial geometries which may be intractable to mesh using hexahedral elements alone. Furthermore, the use of algorithmic approaches such as sum-factorisation reduces the number of floating-point operations, while implementation optimisations such as matrix-free techniques reduce memory bus contention. In introducing simplex element shapes, we have previously investigated, on CPUs, the performance trade-offs of these various algorithmic approaches \cite{Cantwell2011,Cantwell2011h}. The code has recently been ported to GPUs and has been used to explore the impact of various mathematical and algorithmic design choices, such as our treatment of sum-factorisation per output point \cite{Eichstadt2020}, on performance.

Other codes that implement high-order finite element methods have also embraced GPUs, which we acknowledge here in order to give context to our work.
%
% List discretisations and which codes fit under those: acknowledge solid mech (deal.ii and babushka's), but we consider fluid mech. Traditional hexahedral spectral element (patera 1984, Nek5000, NekRS, Neko - hex, lagrange, collocation property, sum fac), mixed-element methods (MFEM, Nektar++)
%
Nek5000 \cite{fischer2007nek5000} is one of the most widely known high-order spectral element codes for the solution of incompressible fluid dynamics. Initial attempts \cite{markidis2015openacc, vincent2022strong} to port this Fortran-based software to heterogeneous computing platforms included the use of the directive-based OpenACC \cite{wienke2012openacc} programming model. However, those early efforts were more recently superseded by NekRS \cite{fischer2022nekrs}, relying on modern C++ and the OCCA \cite{medina2014occa} domain-specific language (DSL), with the support of CUDA, HIP, and OpenCL back-ends. The Neko \cite{jansson2024neko} code is a recent high-order spectral element framework evolved from Nek5000 using modern Fortran and providing direct support for CUDA, HIP and OpenCL back-ends for GPUs, as well as support for general purpose and SX-Aurora vector processors. MFEM \cite{anderson2021mfem} is a general library for leveraging high-order finite-element methods with direct back-end support for CUDA, HIP, and OpenMP, as well as support through OCCA and RAJA \cite{beckingsale2019raja}, while the deal.II \cite{bangerth2007deal} general finite-element library provided initial GPU support \cite{arndt2021deal} using CUDA with more recent efforts \cite{arndt2023deal} through the use of the Kokkos library \cite{kokkos2014}. Additional high-order spectral-element method codes focusing on compressible flow, using either a flux-reconstruction \cite{huynh2007flux} or discontinuous Galerkin \cite{Hesthaven2008, cockburn2000} approaches, are also available. Notably, PyFR \cite{witherden2014pyfr, witherden2025pyfr} relies on a Python front-end combined with a DSL based on Mako templating \cite{bayer2012mako} to provide support for CUDA, HIP, OpenCL, OpenMP, and Metal back-ends, while GAL{\AE}XI \cite{kurz2025galaexi} is a recent extension of the Fortran-based FLEXI \cite {krais2021flexi} discontinuous Galerkin framework providing support for GPUs using CUDA. HORSES3D \cite{ferrer2023high} also provides GPU support with OpenACC \cite{kessasra2024comparison}. Of these codes, only PyFR, deal.II, MFEM, and more recently GAL{\AE}XI/FLEXI provide some support for simplex elements. 

A recent community effort to bring this class of algorithm to GPUs, thereby tackling the scientific need to run such simulations at the exascale, has been investigated through the Centre for Efficient Exascale Discretizations (CEED), a collaboration between some of the above projects, including Nek5000 and MFEM. Through this initiative, a number of benchmark problems have been defined, based on common operations required in finite element solvers.
%\cite{?}% 
These provide a well-defined basis on which to compare different algorithms and implementations. We leverage these benchmarks in this study in order to draw comparisons with previous work.

% TODO: Expand the above to cite other papers which have published on the BK/BP problems.

This paper presents our recent observations concerning the performance of high-order spectral/$hp$ element methods on CPUs and GPUs, for different element shapes and implementation choices. In Section \ref{sec:theory} we discuss the mathematical construction of the methods being considered, including the key finite element operators of interest and expansion bases used within each element. We also document an alternative implementation of the Helmholtz operator which can provide beneficial performance improvements on GPUs and other single instruction, multiple data (SIMD) architectures by reducing the number of inner-product operations. In Section \ref{sec:implementation}, we summarise our generalised approach to managing the different algorithmic implementations within our framework and illustrate how they are efficiently mapped to different back-ends. Section \ref{sec:results} summarises our experiments performed on both CPU and GPU architectures across various element shapes and polynomial degrees and discusses the implications for achieving performance on these platforms. Finally, we conclude with a summary of our observations in Section \ref{sec:summary}.

%% file: Theory.tex
\section{Theory and background}
\label{sec:theory}
In this section, we outline the mathematical framework used in this work by describing a high-order finite element discretisation of the Helmholtz equation. For more details on the method and notation, the reader should consult~\cite{KarniadakisSherwin2005}.
The starting point of our exposition is the weak form of the Helmholtz equation $\nabla^2 u - \lambda u = f$, for an
unknown scalar $u(\boldsymbol{x})$ given a forcing function $f(\boldsymbol{x})$ and the domain $\Omega$.
Using a Galerkin formulation with a condensed (e.g., with boundary conditions enforced) right-hand side this reads as
\begin{equation}
  (\nabla u, \nabla v)_{\Omega} + \lambda (u,v)_{\Omega} = - (f, v)_{\Omega}\,,
\end{equation}
where $(u,v)_\Omega= \int_\Omega u v\dd \Omega$ represents the inner product over the domain (for a vector field, $(\boldsymbol{u},\boldsymbol{v})_\Omega= \int_\Omega \boldsymbol{u} \cdot \boldsymbol{v}\dd \Omega$).  We first outline the basic formulations of spectral elements and then go into the details of the two-dimensional formulation to explain how it affects the computation in practice.

\subsection{Spectral element discretisation}\label{sec:basic}
We adopt an expansion basis of the form $\{\phi_n(\boldsymbol{x}); n=1, \ldots, N_P\}$ on the domain $\Omega$, where each basis function $\phi_n$ represents a polynomial of up to degree $P$. The total number of basis functions, $N_P$, depends on $P$ and element shape, such as $N_P=(P+1)^2$ for the quadrilateral. A function on the element domain is consequently approximated as 
\begin{equation}
  u(\boldsymbol{x}) = \sum_n \hat{u}_n \phi_n(\boldsymbol{x}),
\end{equation}
where $\hat{u}_n$ is the coefficient of each $\phi_n(\boldsymbol{x})$, known as an elemental degree of freedom, while collectively the $\hat{u}_n$ are known as the coefficient space representation of $u(\boldsymbol{x})$. The Galerkin formulation of the Helmholtz equation within an element, $\Omega^e$, reads
\begin{equation}
 \int_{\Omega^e} \sum_{n=1}^{N_P} \left ( \lambda \hat{u}^e_n\phi^e_n(\boldsymbol{x})\phi^e_m(\boldsymbol{x}) + \hat{u}^e_n\nabla\phi^e_n(\boldsymbol{x}) \cdot\nabla\phi^e_m(\boldsymbol{x}) \right ) \dd \Omega
  = -\int_{\Omega^e} f(\boldsymbol{x})\phi^e_m(\boldsymbol{x}) \dd \Omega , \;  m=1, \ldots, N_P \,.
\end{equation}
The superscript $e$ represents individual element, but for brevity, we will omit it in the basis notation $\phi_n$ in the later context.  We note that at this stage, we have not associated $\phi_n(x)$ with any particular set of collocating points; in other words, they are general-type functions and may be nodal, modal or even orthogonal \cite{karniadakis-2013}. Given a set of unique points in an element, we can build a Lagrange expansion that spans the same space to transform between modal and nodal representations. We will exploit this later in this work.

Since the coefficients $\hat{u}^e_n$ are independent of $\boldsymbol{x}$, we can factor them out and rewrite the equations as a matrix system
\begin{equation}
 \bm{H}^e \, {\boldsymbol{ \hat{ {u}} }} = {\boldsymbol{ \hat{ {f}} }}\,,
\end{equation}
with the elemental Helmholtz operator $\bm{H}^e$, the vector of solution coefficients ${\boldsymbol{ \hat{ {u}} }}$, and the right-hand side vector $\hat{ \boldsymbol{f}}$. As a matter of practice, we define the polynomial basis on the reference element $\Omega_{\text{st}}$ (as seen in the next section), so the inner product will be performed on the reference element with the knowledge of the one-to-one mapping $\boldsymbol{x} = \chi^e(\boldsymbol{\xi})$ between global Cartesian coordinates $\boldsymbol{x}$ to local Cartesian coordinate $\boldsymbol{\xi} $. Then we obtain
\begin{equation}
\bm{H}^e[m][n] = \int_{\Omega_{\text{st}}} \left[\lambda\phi_n(\boldsymbol{\xi})\phi_m(\boldsymbol{\xi}) +
    (\bm{J}^{e})^{-T}\nabla\phi_n(\boldsymbol{\xi}) \cdot
    (\bm{J}^{e})^{-1}\nabla\phi_m(\boldsymbol{\xi})\right] |\bm{J}^e|\, \dd \Omega\,.
\end{equation}
Here we have applied the chain rule
\begin{equation}
\nabla\phi( \boldsymbol{x}) = \frac{\partial \boldsymbol{\xi}}{\partial \boldsymbol{x}}\frac{\partial \phi}{\partial  \boldsymbol{\xi}} = \left( \bm{J}^e \right)^{-1}\nabla\phi(\boldsymbol{\xi}) \,,
\end{equation}
where $\bm{J}^e = \nabla \chi^e(\boldsymbol{\xi})$ is the Jacobian matrix of the mapping, which can be constant within an element if the mapping is {\em affine} or a polynomial of degree $N$ when accomplishing iso-parametric analysis. Separating the operator implementations for these two different types of mappings is crucial for achieving maximal performance.

The integral in the above formulations is evaluated using Gaussian quadrature.  With weights $w_l$ and quadrature points $\boldsymbol{\xi}_l\in\Omega_{\text{st}}$ we can write the discrete elemental Helmholtz operator as
\begin{equation}
\bm{H}^e_\delta[m][n]  = \sum_ l  \left[\lambda w_l \phi_n(\boldsymbol{\xi}_l)\phi_m(\boldsymbol{\xi}_l) +
  w_l  (\bm{J}^{e}_l)^{-T}\nabla\phi_n(\boldsymbol{\xi}_l) \cdot
    (\bm{J}^{e}_l)^{-1}\nabla\phi_m(\boldsymbol{\xi}_l)\right] |\bm{J}^e_l|\,.
\end{equation}
The multi-dimensional quadrature is constructed using a tensor-product rule. We typically choose at least $P+2$  Gauss-Lobatto points along each coordinate for expansions of order $P$ in order to be able to exactly evaluate the product of two polynomial shape functions $\phi_p(\xi)\phi_q(\xi)$ in the standard region.

Using the matrix notation of \cite{KarniadakisSherwin2005}, the discrete Helmholtz operator can be further split into the mass operator $\bm{M}^e_\delta$ and the Laplacian operator $\bm{L}^e_\delta$ yielding
\begin{equation}
\bm{H}^e_\delta  = \lambda \bm{M}^e_\delta + \bm{L}^e_\delta  ~.
\end{equation}

\subsection{Two-dimensional formulation} \label{sec:detail}
For two-dimensional element quadrilaterals and triangles, the local Cartesian coordinates are denoted as  $\boldsymbol{\xi}=(\xi_1,\xi_2)$. For quadrilateral regions we use a tensor-product expansion basis of the form 
$\{\phi_{pq}(\xi_{1},\xi_{2})\ ;\ p=0, \ldots, P\ ;\ q=0, \ldots,  P \}$. 
Triangular elements will be constructed using the tensor-product 
approximations as described in Section~\ref{sec:hk_method_tri}. The approximated function is then given by the summation
\begin{equation}
  u(\bm{\xi}) = \sum_{m(p,q)} \hat{u}_{m} \phi_{m}(\xi_1,\xi_2) = \sum_{p=0}^P \sum_{q=0}^P \hat{u}_{pq} \phi_{pq}(\xi_{1},\xi_{2}),
\end{equation}
where the indices $m(p,q)$ are obtained as a linear combination of indices of $p,q$ depending on the type of element considered in the approximation. The two-dimensional discrete elemental Helmholtz operator now reads
\begin{equation}
\begin{split}
\bm{H}^e_\delta[m][n]  
&= \lambda  \sum_i \sum_j  w_i w_j \phi_{pq}(\xi_{1i},\xi_{2j})\phi_{rs}(\xi_{1i},\xi_{2j}) |\bm{J}^e_{ij}|  \\
&+         \sum_i \sum_j w_i w_j  (\bm{J}^{e}_{ij})^{-T}\nabla\phi_{pq}(\xi_{1i},\xi_{2j}) \cdot
    (\bm{J}^{e}_{ij})^{-1}\nabla\phi_{rs}(\xi_{1i},\xi_{2j}) |\bm{J}^e_{ij}|  \,,
\end{split}
    \label{eq:helm}
\end{equation}
where the indices $i$ and $j$ again span over the quadrature points, and the index $n(r,s)$ is a linear combinations of the indices  $r$, and $s$. For specific details of this calculation, the interested reader 
could consult~\cite{KarniadakisSherwin2005} (Section 4.1.5).

The mass operator can now be rewritten in matrix form as
\begin{equation}\label{eq:mass}
	\bm{M}^e_\delta = \bm{B}^T \bm{W} \bm{B} \,, 
\end{equation}
where the dense basis matrix $\bm{B}$ contains terms of the form $\phi_{pq}(\xi_{1i},\xi_{2j})$, and 
the components of the diagonal matrix ${\bm{W}}$ are the weights and the Jacobian $w_i w_j |\bm{J}^e_{ij}|$.  Two comments are in order: firstly, the matrix $\bm{B}$ can be considered a basis evaluation matrix which yields a vector containing the expansion evaluated at the points $(\xi_{1i},\xi_{2j})$ when acting upon a vector containing the expansion coefficients ${\boldsymbol{\hat{u}}}$. As presented above, these points denote the quadrature points for evaluating the inner product. In nodal spectral elements, these (quadrature) points can also denote the positions at which a collocating Lagrange basis can be built when one restricts oneself to using the same points for the nodal expansion and the quadrature.  It is then important to appreciate that the matrix $\bm{B}$ can alternatively be viewed as a transformation from a modal, or a different nodal, basis to a nodal basis at the quadrature points over the local coordinate system constructed at the quadrature points. We will exploit this observation below to efficiently organise our computation.

We now present three necessary building blocks.  Firstly, the derivative of the basis functions with respect to the local coordinates for two-dimensional expansions is
\begin{equation}
	\nabla\phi(\xi_{1i},\xi_{2j}) = \left(\begin{matrix}
    \dfrac{\partial \phi}{\partial \xi_{1}}  \\[2ex]
    \dfrac{\partial \phi}{\partial \xi_{2}}
  \end{matrix}\right)_{(\xi_{1i},\xi_{2j})}.
\end{equation}

\noindent Secondly, the inverse of the Jacobian matrix $(\bm{J}^e_{ij})^{-1}$ is given by:
\begin{equation}
	(\bm{J}^e_{ij})^{-1} = \left(\begin{matrix}
    \dfrac{\partial \xi_1}{\partial x_1} & \dfrac{\partial \xi_2}{\partial x_1} \\[2ex]
    \dfrac{\partial \xi_1}{\partial x_2} & \dfrac{\partial \xi_2}{\partial x_2}
  \end{matrix} \right)_{(\xi_{1i},\xi_{2j})}.
\end{equation}

As noted in \cite{KarniadakisSherwin2005}, when the elemental mapping function is not affine, the entries of this matrix represent non-polynomial, yet smooth functions. Higher quadrature order might be required to see spectral convergence when these terms deviate significantly from being polynomial \cite{mengaldo-2015}. The Laplacian operator can now be expanded as
\begin{equation}
\begin{split}
\bm{L}^e_\delta = \sum_i \sum_j w_i w_j  \Biggl[ 
\biggl( \frac{\partial \xi_1}{\partial x_1} \frac{\partial \phi_{pq}}{\partial \xi_1} 
&+      \frac{\partial \xi_2}{\partial x_1} \frac{\partial \phi_{pq}}{\partial \xi_2} \biggr) 
\biggl( \frac{\partial \xi_1}{\partial x_1} \frac{\partial \phi_{rs}}{\partial \xi_1} 
+       \frac{\partial \xi_2}{\partial x_1} \frac{\partial \phi_{rs}}{\partial \xi_2} \biggr) \\ +
\biggl( \frac{\partial \xi_1}{\partial x_2} \frac{\partial \phi_{pq}}{\partial \xi_1} 
&+      \frac{\partial \xi_2}{\partial x_2} \frac{\partial \phi_{pq}}{\partial \xi_2} \biggr) 
\biggl( \frac{\partial \xi_1}{\partial x_2} \frac{\partial \phi_{rs}}{\partial \xi_1} 
+       \frac{\partial \xi_2}{\partial x_2} \frac{\partial \phi_{rs}}{\partial \xi_2} \biggr) \Biggr]_{ij} |\bm{J}^e_{ij}| \,.
\end{split}
\end{equation}
Recall the basis matrix relation mentioned earlier. Any polynomial basis function can be represented as a Lagrange expansion $h_{m(k,l)}(\xi_1,\xi_2) = h_k(\xi_1)h_l(\xi_2)$ of appropriate order based on the function values at the quadrature points:
\begin{equation}\phi(\xi_{1},\xi_{2})= \sum_{m(k,l)}  h_m(\xi_1,\xi_2) \phi(\xi_{1k}, \xi_{2l}) = \sum_{k}\sum_{l} h_k(\xi_1)h_l(\xi_2)\phi(\xi_{1k},\xi_{2l}) \,,
\label{eq:phitoh0}
\end{equation}
where $h_k(\xi_1)$ is the one-dimensional Lagrange polynomial corresponding to the nodal point $\xi_{1k}$ and similarly for $h_l(\xi_2)$. The multi-dimensional Lagrange polynomial can always be constructed by the tensor product, regardless of element shape. Then we can rewrite the derivative of the basis functions into
\begin{equation}
     \dfrac{\partial \phi}{\partial \xi_{1}}\left(\xi_{1},\xi_{2}\right) = \sum_{m(k,l)}  \dfrac{\dd h_m(\xi_1,\xi_2)}{\dd\xi_1} \phi(\xi_{1k}, \xi_{2l}) = \sum_{k}\sum_{l} \frac{\dd h_k(\xi_{1})}{\dd \xi_1} h_l(\xi_2) \phi(\xi_{1k},\xi_{2l}) \,.
     \label{eq:dphitodh}
\end{equation}
\noindent The derivatives of functions at quadrature points can now be expressed as
\begin{equation}
     \Biggl[\dfrac{\partial \phi_{pq}}{\partial \xi_{1}}\Biggr]_{ij}=  \sum_{k}\frac{\dd h_k(\xi_{1i})}{\dd \xi_{1}} \phi_{pq}(\xi_{1k},\xi_{2j})=\bm{D}_{\xi_1}\bm{B} \,,
\end{equation}
where we have applied the collocation property $h_l(\xi_{2j})=\delta_{lj}$, because there is an implicit assumption (made explicit here) that the points at which the basis matrix is evaluated are the same as the points at which the differentiation matrix is formed. $\bm{D}_{\xi_1}$ is actually a block-diagonal differentiation matrix with one-dimensional derivatives of Lagrange polynomials $\frac{\dd h_k(\xi_{1i})}{\dd \xi_{1}}$ as entries, known as {\em collocation differentiation} \cite{karniadakis-2013} (p.65). The action of ${\bm {D B}}$ on an expansion vector ${\boldsymbol{\hat{u}}}$ can be interpreted (reading right to left) as the basis matrix transforming the expansion coefficients to the corresponding nodal coefficients of a Lagrange expansion, followed by the differentiation matrix action on the function at the nodal values.

As acknowledged in \cite{KarniadakisSherwin2005} (p.178), through substitutions and rearrangement, you arrive at the following expression:

\begin{eqnarray}  
{\bm L}^e &=& \bm{B}^T \biggl[{\bm  D}^T_{\xi_1}\bm{\Lambda} \left( \frac{\partial \xi_1}{\partial x_1} \right) 
+      {\bm D}^T_{\xi_2}         \bm{\Lambda} \left( \frac{\partial \xi_2}{\partial x_1} \right) \biggr]\,{\bm W}\, \biggl[\bm{\Lambda} \left( \frac{\partial \xi_1}{\partial x_1} \right){\bm D_{\xi_1}} 
+     \bm{\Lambda} \left( \frac{\partial \xi_2}{\partial x_1} \right) {\bm D}_{\xi_2} \biggr] {\bm B} \nonumber\\
&+& {\bm B}^T\biggr[{\bm D}^T_{\xi_1}\bm{\Lambda} \left( \frac{\partial \xi_1}{\partial x_2} \right)  
+               {\bm D}^T_{\xi_2}\bm{\Lambda} \left( \frac{\partial \xi_2}{\partial x_2} \right)\biggr] \,{\bm W}\,\biggr[\bm{\Lambda} \left( \frac{\partial \xi_1}{\partial x_2} \right)  {\bm D}_{\xi_1}
+               \bm{\Lambda} \left( \frac{\partial \xi_2}{\partial x_2} \right){\bm D}_{\xi_2}\biggr]\bm{B} \,, \label{eq:Lap}
\end{eqnarray}

\noindent where $\bm{\Lambda}(f)$ represents a diagonal matrix with entries  $f(\xi_{1i},\xi_{2j})$. This formulation also demonstrates that the symmetry of the local operator is maintained.

\subsubsection{Alternate interpretations of the Helmholtz operator}\label{sec:reinterp}

%As noted earlier ${\bm W}$ in Equation (\ref{eq:Lap}) is also a diagonal matrix and hence can be commuted to the left, outside the bracketed notation. 
Combining Equation (\ref{eq:Lap}) and Equation (\ref{eq:mass}), we observe that we can factor out $\bm{B}^T$ and $\bm{B}$ to obtain the discrete Helmholtz matrix:
\begin{eqnarray}  
\bm{H}^e_\delta = \bm{L}^e_\delta  + \lambda \bm{M}^e_\delta &=& \bm{B}^T\left \{ \biggl[{\bm D}^T_{\xi_1}\bm{\Lambda} \left( \frac{\partial \xi_1}{\partial x_1} \right) 
+      {\bm D}^T_{\xi_2}         \bm{\Lambda} \left( \frac{\partial \xi_2}{\partial x_1} \right) \biggr]\,  {\bm W} \,\biggl[\bm{\Lambda} \left( \frac{\partial \xi_1}{\partial x_1} \right){\bm D_{\xi_1}} 
+     \bm{\Lambda} \left( \frac{\partial \xi_2}{\partial x_1} \right) {\bm D}_{\xi_2} \biggr] \right . \nonumber\\
&+&  \left .  \biggr[{\bm D}^T_{\xi_1}\bm{\Lambda} \left( \frac{\partial \xi_1}{\partial x_2} \right)  
+   {\bm D}^T_{\xi_2}\bm{\Lambda} \left( \frac{\partial \xi_2}{\partial x_2} \right)\biggr] \,{\bm W}\,\biggr[\bm{\Lambda} \left( \frac{\partial \xi_1}{\partial x_2} \right)  {\bm D}_{\xi_1}
+               \bm{\Lambda} \left( \frac{\partial \xi_2}{\partial x_2} \right){\bm D}_{\xi_2}\biggr]+\lambda\bm{W} \right \} \bm{B} \,. \label{eq:helm1}
\end{eqnarray}

We observe that the matrix operator ${\bm B^T}$
is similar to the inner product with respect to the basis $\phi_{pq}(\xi_1,\xi_2)$ operation, although the full operations also requires a diagonal matrix ${\bm W}$. However, in Equation (\ref{eq:helm1}), there is only one ``inner product" type operation rather than two in Equation (\ref{eq:Lap}),  one with respect to 
$\frac{\partial \phi_{pq}}{\partial \xi_{1} }$ 
(related to the matrix $\bm{B}^T\bm{D}^T_{\\\xi_1}$) and the other with respect to $\frac{\partial \phi_{pq}}{\partial \xi_2}$ (related to the matrix $\bm{B}^T\bm{D}^T_{\\\xi_2}$). This seemingly rather trivial observation has an important implication on SIMD architectures since the inner product type operation can be quite costly, particularly for mixed shape expansions discussed in the next section. Even for quadrilateral and hexahedral modal expansions (or nodal expansions of a different polynomial degree than the order of the quadrature that is adopted), it means we can perform the operation within the curly brackets in Equation (\ref{eq:helm1}) as if we were using a nodal expansion with collocation properties but pre- and post-multiply by a  $\bm{B}^T$ and $\bm{B}$
matrix respectively. Indeed, in our implementation, basing our formulation on Equation (\ref{eq:helm1})  led to a percentage improvement of between 20\% to 50\%, depending on polynomial order, on hexahedral elements and 50\% on tetrahedral elements (see Section \ref{sec:coll}).

Consideration of Equation (\ref{eq:helm1}) highlights that  the  Helmholtz operator can be decomposed into a transformation from a modal basis  (or indeed any other basis not at the quadrature points) to a tensor product Lagrange basis defined at a set of Gaussian quadrature points (represented by the matrix $\bm{B}$), followed by a series of collocated or even point-wise operations at the quadrature points,  and finally a ``restriction" back to the modal basis (represented by the matrix $\bm{B}^T$). This is exact, as long as the span of the Lagrange polynomial basis is no less than the polynomial space or the original expansion $\phi_n(\bf \bf \xi)$, which is almost always true since it is common practice to adopt a quadrature order that is never lower than the basis order. This formulation can also be applied to any operator that is tested against the expansion basis or its derivatives, such as the linear advection-diffusion-reaction operator. Although it is quite straightforward to see what is happening at the discrete level, it is not immediately so obvious what is being evaluated at the analytical mathematical level initially outlined in Section \ref{sec:basic}.

\subsubsection*{Trasnsformation versus evaluation of the Basis}
We observe that the terms in the curly brackets of Equation (\ref{eq:helm1}) is the discrete form of the weak Laplacian of the Lagrange basis through the quadrature points, i.e. $(\nabla h_{n}, \nabla h_{m})$. This once again leads us to observe that the discrete (or continuous) weak Laplacian with respect to the expansion basis can be constructed as the weak Laplacian with respect to any Lagrange basis with an equivalent or larger span, pre- and post-multiplied by the associated expansion matrices $\bm{B}^T$ and $\bm{B}$.

To show this relationship more explicitly, we once again start by representing our polynomial $\phi_{m}(\xi_1,\xi_2)$ in terms of a Lagrange tensor product expansion $h_{n(i,j)}(\xi_1,\xi_2) = h_i(\xi_{1}) h_j(\xi_{2})$ of larger span. As already shown in Equation (\ref{eq:phitoh0}), we can identify a set of coefficients $b_{n}$ that relate the expansion basis, $\phi_{m}$ to the Lagrange expansion, so that 
\begin{equation}
  \phi_{m}(\xi_1,\xi_2)  = \sum_{n} h_{n}(\xi_1,\xi_2)\, b_n.
  \label{eq:phitoh}
\end{equation}

For the specific case we have considered up to now, the Lagrange polynomial is defined through a set of quadrature points. Applying the collocation property of the  Lagrange polynomial, it is evident that the entries of the columns of the matrix $\bm{B}$ are the same as $b_{n}$, i.e. $\bm{B}[n,m] = b_n$. However, we note that in Equation (\ref{eq:phitoh}), $b_{n}$ can be interpreted as the transformation of the  expansion coefficients relating the Lagrange basis to $\phi_{m}$, instead of the basis values evaluated at the quadrature points. 

We recall that the Laplacian operator requires an inner product with respect to the local derivative of the basis $(\nabla \phi_{m},\nabla \phi_{n})$. Taking the gradient of Equation (\ref{eq:phitoh}) we  obtain
 \begin{equation}
  \nabla \phi_{m}(\xi_1,\xi_2)  = \sum_{n} \nabla h_{n}(\xi_1,\xi_2)\, b_n= \sum_{n} \nabla h_{n}(\xi_1,\xi_2)\bm{B}[n,m]\,.
\label{eq:GradPhitoh}
\end{equation}
Finally, each component of the weak Laplacian matrix can then be expressed as
\begin{equation}
\bm{L}[m,n] = (\nabla \phi_m, \nabla \phi_n) = \left (\sum_{r} \bm{B}[r,m] \nabla h_{r},\sum_{s} \bm{B}[s,n] \nabla h_{s}\right ) = (\bm{B}[:,m])^T \bm{L}^h \bm{B}[:,n]
\label{eq:helm2}
\end{equation}
where
\[
\bm{L}^h[m,n] = \left (\nabla h_m, \nabla h_n \right),
\]
which leads us to Equation (\ref{eq:helm1}) when discretised and expressed in matrix form, but highlights the relationship between the continuous components of $\bm{L}$ and $\bm{L}^h$. In Equation (\ref{eq:helm1}) we interpret $\bm{B}^T$ as an evaluation of the basis $\phi_m$ at the quadrature points; however in Equation (\ref{eq:helm2}) we can interpret the role of the $\bm{B}$ matrix as a transformation from the $\phi_m$ expansion to the Lagrange basis, $h_n$, through the the quadrature points.

 %We can therefore appreciate that the building blocks of the Laplacian operator depend on the evaluation of a local inner product of the form: 
%\begin{equation}
%    \left ( \frac{\partial \phi_{pq}}{\partial \xi_1}, f \right )_{\Omega_e}. 
%\label{eq:iprod}
%\end{equation}
% Finally, we observe that we can recast Equation (\ref{eq:iprod}) as 
%\[
%    \left ( \frac{\partial \phi_{pq}}{\partial \xi_1}, f \right )_{\Omega_e} = 
%\sum_{ij}b_{ij}\left ( \frac{\partial h_{ij}}{\partial \xi_1}, f \right )_{\Omega_e} = \sum_{n(i,j)} \bm{B}^T[n][m] \left( \frac{\partial h_{ij}}{\partial \xi_1}, f \right )_{\Omega_e} \,,
%\]
%where as used previously $n(i,j)$ and $m(p,q)$ are the indexing relating $i,j$ and $p,q$ to the entries of the $\bm{B}$ matrix. 

%This therefore also leads us to the interpretation that Equation (\ref{eq:Lap1}) can be understood as an inner product with respect to the derivative of the Lagrange basis (which is encapsulated in the matrices $\bm{D}_{\xi_1}, \bm{D}_{\xi_2}$), which is then transformed into the inner product with respect to the derivative of the expansion basis $\phi_{pq}$ by the matrix $\bm{B}^T$. Interestingly, the role of the diagonal matrix ${\bm W}$ , strictly speaking,  remains with the inner operators of the Lagrange basis in the analytic interpretation, although it was the commutation of this diagonal matrix at the discrete level that triggered our exploration of this reinterpretation.

\subsection{Extension to hybrid elements for unstructured meshes}
\label{sec:hk_method_tri}

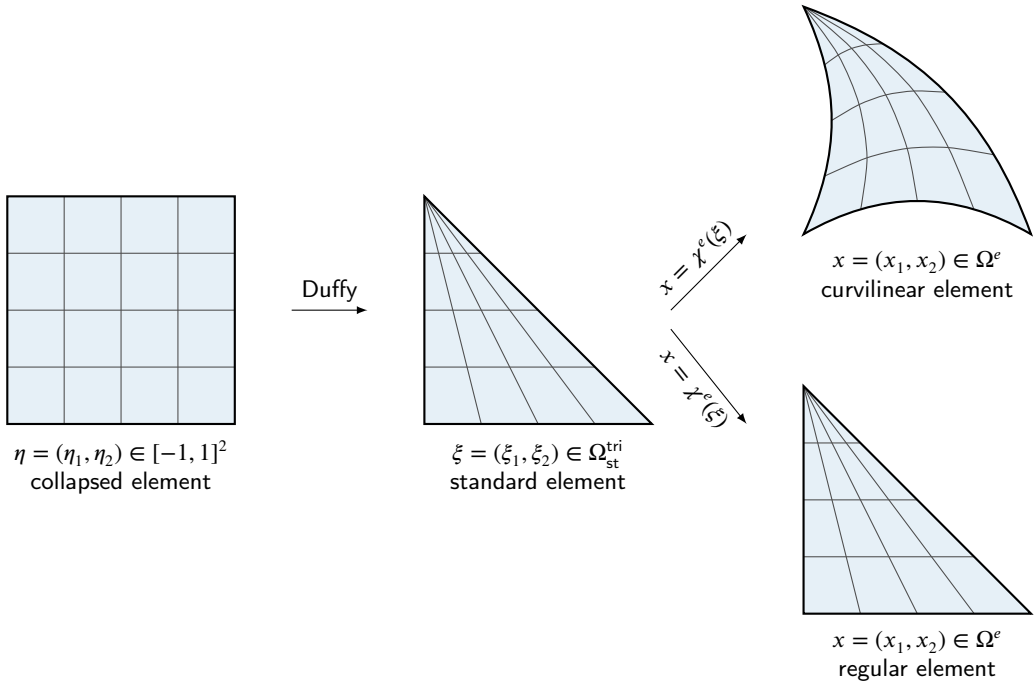
\begin{figure}[tb!]
  \begin{center}
    \begin{tikzpicture}
   	  % collapsed element
	  \begin{scope}[xshift=-1.5cm]	      
	      \fill[fill=blue!10] (0,0) rectangle +(3,3);
	      \draw[black!70] (0.75,0) -- (0.75,3);
	      \draw[black!70] (1.50,0) -- (1.50,3);
	      \draw[black!70] (2.25,0) -- (2.25,3);
	      \draw[black!70] (0,0.75) -- (3,0.75);
	      \draw[black!70] (0,1.50) -- (3,1.50);
	      \draw[black!70] (0,2.25) -- (3,2.25);
	      \draw[thick] (0,0) rectangle +(3,3);
	      \node[align=center,below] at (1.5,-0.1) {$\eta = (\eta_1,\eta_2) \in [-1,1]^2$\\collapsed element};
      \end{scope}
      \draw[-latex] (2.25,1.5) to node[above] {Duffy} (3.25,1.5);
      
      %standard element
      \fill[fill=blue!10] (4,0) -- (7,0) -- (4,3) -- cycle;
      \draw[black!70] (4.75,0) -- (4,3);
      \draw[black!70] (5.50,0) -- (4,3);
      \draw[black!70] (6.25,0) -- (4,3);
      \draw[black!70] (4,0.75) -- (6.25,0.75);
      \draw[black!70] (4,1.50) -- (5.50,1.50);
      \draw[black!70] (4,2.25) -- (4.75,2.25);
      \draw[thick] (4,0) -- (7,0) -- (4,3) -- cycle;
      \node[align=center,below] at (5.5,-0.1) {$\xi = (\xi_1,\xi_2) \in \Omega^{\text{tri}}_{\text{st}}$\\standard element};

      \draw[-latex] (7.25,1.5) to node[above,midway,sloped] {$x = \chi^e(\xi)$} (8.25,2.5);
      \draw[-latex] (7.25,1.25) to node[below,midway,sloped] {$x = \chi^e(\xi)$} (8.25,0);

      % curvilinear element
      \begin{scope}[xshift=9cm, yshift=2.5cm]
        \node[align=center,below] at (1.5,-0.1) {$x = (x_1,x_2) \in \Omega^e$\\curvilinear element};
        \fill[fill=blue!10]
        (0,0) 
        to [bend left] coordinate[pos=0.25] (A1) coordinate[pos=0.5] (A2) coordinate[pos=0.75] (A3) (3,0)
        to [bend right=25] coordinate[pos=0.25] (B1) coordinate[pos=0.50] (B2) coordinate[pos=0.75] (B3) (0,3)
        to [bend left=25] coordinate[pos=0.25] (C1) coordinate[pos=0.50] (C2) coordinate[pos=0.75] (C3) cycle;

        \path (B1) to [bend right=20] coordinate[pos=0.25] (D1) coordinate[pos=0.5] (D2) coordinate[pos=0.75] (D3) (C3);
        \path (B2) to [bend right=20] coordinate[pos=0.25] (E1) coordinate[pos=0.5] (E2) coordinate[pos=0.75] (E3) (C2);
        \path (B3) to [bend right=20] coordinate[pos=0.25] (F1) coordinate[pos=0.5] (F2) coordinate[pos=0.75] (F3) (C1);

        \draw[black!70] plot[smooth] coordinates {(A1) (D3) (E3) (F3) (0,3)};
        \draw[black!70] plot[smooth] coordinates {(A2) (D2) (E2) (F2) (0,3)};
        \draw[black!70] plot[smooth] coordinates {(A3) (D1) (E1) (F1) (0,3)};
        \draw[black!70] plot[smooth] coordinates {(B1) (D1) (D2) (D3) (C3)};
        \draw[black!70] plot[smooth] coordinates {(B2) (E1) (E2) (E3) (C2)};
        \draw[black!70] plot[smooth] coordinates {(B3) (F1) (F2) (F3) (C1)};

        \draw[thick] (0,0) to [bend left] (3,0) to [bend right=25] (0,3) to [bend left=25] cycle;
      \end{scope}

      % regular element
      \begin{scope}[xshift=5cm,yshift=-2.5cm]
        \node[align=center,below] at (5.5,-0.1) {$x = (x_1,x_2) \in \Omega^e$\\regular element};
      \fill[fill=blue!10] (4,0) -- (7,0) -- (4,3) -- cycle;
      \draw[black!70] (4.75,0) -- (4,3);
      \draw[black!70] (5.50,0) -- (4,3);
      \draw[black!70] (6.25,0) -- (4,3);
      \draw[black!70] (4,0.75) -- (6.25,0.75);
      \draw[black!70] (4,1.50) -- (5.50,1.50);
      \draw[black!70] (4,2.25) -- (4.75,2.25);
      \draw[thick] (4,0) -- (7,0) -- (4,3) -- cycle;

      \end{scope}
    \end{tikzpicture}
  \end{center}
  \caption{ Representation of a high-order triangular element with three different coordinate systems. Elements may be curvilinear or benefit from reduced geometric information if regular straight-sided.
The distribution of quadrature points is illustrated with equispaced lines.}
  \label{fig:collapsed}
\end{figure}

For triangles in 2D and tetrahedra, pyramids, and prisms in 3D, to obtain a tensor-product construction of our local elemental operations, 
we employ the standard square-to-triangle Duffy transformation~\cite{Duffy1982} to define two independent coordinate directions. In the case of triangles, only a single use of this transformation is needed. For the 3D shapes, multiple successive actions of the Duffy transformation are required as outlined in ~\cite{KarniadakisSherwin2005}. We will limit ourselves to describing things on a triangle as the other shapes follow quite naturally.

Using a Duffy transformation allows us to define a basis on a tensor-product square with coordinates $\eta\in[-1,1]^2$ and map this to a triangle with reference coordinates $\xi\in\Omega_{\text{st}} = \{ (\xi_1, \xi_2)\ ;\ \xi_{1,2}\in[-1,1],\ \xi_1+\xi_2 \leq 0 \}$.
This transformation is illustrated in Figure~\ref{fig:collapsed}. As previously mentioned, whereas tensor-product quadrature operations were previously accomplished in the standard element for quadrilaterals, tensor-product operations will now be accomplished on the collapsed coordinate system. Analytically, the Duffy mapping is defined as
\begin{equation}
  \eta_1 = 2\frac{1+\xi_1}{1-\xi_2}-1, \quad \eta_2 = \xi_2 \,.
\end{equation}
For triangles, Equation~\eqref{eq:helm} is solved on the collapsed coordinate space, instead of the Cartesian coordinate space. 
We hence require a partial derivative with respect to the Cartesian space, which is determined by the chain rule as
\begin{equation}
  \left(\begin{matrix}
    \dfrac{\partial}{\partial\xi_1}  \\[3ex]
    \dfrac{\partial}{\partial\xi_2}
  \end{matrix} \right) = \left(\begin{matrix}
    \dfrac{2}{1-\eta_2} & 0 \\[3ex]
    2\dfrac{1+\eta_1}{1-\eta_2} & 1
  \end{matrix} \right) \left(\begin{matrix}
    \dfrac{\partial}{\partial\eta_1}  \\[3ex]
    \dfrac{\partial}{\partial\eta_2}
  \end{matrix} \right)\, ,
\end{equation}
or in matrix notation as
\begin{equation}
  \nabla_{\xi}  = \bm{G} \,\nabla_{\eta}\,.
\end{equation}

As discussed in \cite{KarniadakisSherwin2005} when designing tensor products expansions for simplex shapes  we use a more generalised form. So for triangles the generalised product is of the form: 
\[
\phi^{\text{tri}}_{pq}(\eta_i,\eta_2) = \psi^a_{q}(\eta_1) \psi^b_{pq}(\eta_2)
\]
where $0 \leq p \leq P$ and $0 \leq p+q \leq P$. The dependency of the second product on both $p$ and $q$ as well as the restriction that $p+q \leq P$ makes the inner product operator more challenging to implement on a SIMD environment. The interpretation provided in section \ref{sec:reinterp} is therefore even more relevant when we consider complex operators such as the Laplacian operator in hybrid elements. By using a Lagrange tensor product basis $h_p(\eta_1)h_q(\eta_2)$, similar to the quadrilateral case,  we restricted the cumbersome triangular space  $\phi^{\text{tri}}_{pq}$ within the matrix operation $\bm{B}^T$.

%% file: Implementation.tex
\section{Implementation}
\label{sec:implementation}
In this section, we outline key design considerations that enable spectral/$hp$ element methods to be made performant across host and accelerators. 
% %
Through our investigations, three key concepts have been identified for the spectral/$hp$ element operations in order to deliver a robust and performant implementation: (i) use of a memory space model, (ii) use of an execution space model, and (iii) support for multiple implementation strategies per operator. In the following section, we outline these concepts, focusing specifically on the implementation strategies considered in the results presented in Section~\ref{sec:results}. Implementation details, in the context of Nektar++, are given in Appendix~\ref{app:algo}.

\subsection{Memory space model}
\label{sec:implementation:mem}
As high-performance computing systems have moved to architectures involving accelerator devices with discrete memory, appropriate memory management between the host CPU and the device has become an important factor impacting performance. 
% %
The concept of a {\it memory space} represents the physical location of the memory in which the data reside. We note that a similar memory space model is also adopted by the \emph{Kokkos} library \cite{edwards2014kokkos}. Current computer architectures necessitate two memory spaces: 
\begin{itemize}
\item \textbf{HostSpace}: Represents the memory accessed from the host, typically the CPU, including those with vector units.
\item \textbf{DeviceSpace}: Represents the memory accessed from an accelerator device, typically a GPU with discrete memory. 
\end{itemize}
%%
%Figure \ref{fig:memory_space} provides a schematic representation of the memory space model we have adopted where the double-ended arrow represents data movements. 

%As detailed below, a semi-automatic memory management system has been developed to determine when host/device memory synchronisation is required without explicit function calls for memory movements.
%%
%\begin{figure}
%\centering
%\includegraphics[width=0.6\textwid%th]{images/memory_space.pdf}
%\caption{Memory space model %highlighting a semi-automatic %memory management between the %local Device (i.e. GPU) and Host %%(i.e. CPU)}
%\label{fig:memory_space}
%\end{figure}
%%
%
In a heterogeneous computing framework, the memory movement and synchronisation between a host and a device usually presents a significant challenge to the programmer. Incorporating memory access qualifiers enables the automated transfer of data between host and device as necessary, relieving the developer of much of this burden. In particular, the following types of access should be supported:
\begin{itemize}
	\item \textbf{ReadOnly}: Requires that the memory has already been initialised in the requested memory space. Data will only be read from memory. Memory transfer occurs if the data in the current memory space was marked as invalid.
	\item \textbf{WriteOnly}: Memory is allocated if needed and marked as initialised. No memory transfer occurs. The data in the requested memory space is now marked as valid, while the data in the other memory space is marked as invalid.
	\item \textbf{ReadWrite}: Requires that the memory has already been initialised in the requested memory space. Memory transfer occurs if the data on the requested memory space was marked as invalid. The data in the requested memory space is then marked as valid, while the data in the other memory space is marked as invalid.
\end{itemize}

Similar approaches based on a memory access qualifier model have been used by parallel programming frameworks such as {\it OpenCL} and {\it SYCL} or by parallel programming libraries such as {\it MFEM} \cite{anderson2021mfem} and {\it StarPU} \cite{augonnet2009starpu}. More details on the implementation of memory access specifiers are given in Appendix~\ref{app:algo}.

\subsection{Execution space model}
\begin{figure}[htb]
\centering
\begin{forest}
  for tree={
      font=\sffamily\bfseries,
      line width=1pt,
      draw=linecol,
      ellip,
      align=center,
      child anchor=north,
      parent anchor=south,
      drop shadow,
      l sep+=12.5pt,
      edge path={
        \noexpand\path[color=linecol, rounded corners=5pt, >={Stealth[length=10pt]}, line width=1pt, ->, \forestoption{edge}]
          (!u.parent anchor) -- +(0,-5pt) -|
          (.child anchor)\forestoption{edge label};
        },
  }
  [Execution Space, inner color=col1in, outer color=col1out, calign=child, calign child=2
    [Serial, inner color=col3in, outer color=col3out]
    [AVX, inner color=col2in, outer color=col2out, calign=child, calign child=3
      [SSE2, inner color=col3in, outer color=col3out, edge=dashed]
      [AVX2, inner color=col3in, outer color=col3out, edge=dashed]
      [AVX512, inner color=col3in, outer color=col3out, edge=dashed]
      [SVE/SVE2, inner color=col3in, outer color=col3out, edge=dashed]
    ]
    [Device, inner color=col2in, outer color=col2out, calign=child, calign child=2
      [CUDA, inner color=col3in, outer color=col3out, edge=dashed]
      [HIP, inner color=col3in, outer color=col3out, edge=dashed]
      [SYCL, inner color=col3in, outer color=col3out, edge=dashed, calign=child, calign child=2
       [CUDA, inner color=col5in, outer color=col5out, edge=dashed]
       [HIP, inner color=col5in, outer color=col5out, edge=dashed]
       [OpenCL, inner color=col5in, outer color=col5out, edge=dashed]
       ]
      ]
    ]
  ]
\end{forest}
\caption{Nektar++ execution space model hierarchy. Execution spaces are shown in green, while specific backends are shown in blue. Availability of backends depends on the host and/or device architecture.}
\label{fig:execution_space}
\end{figure}
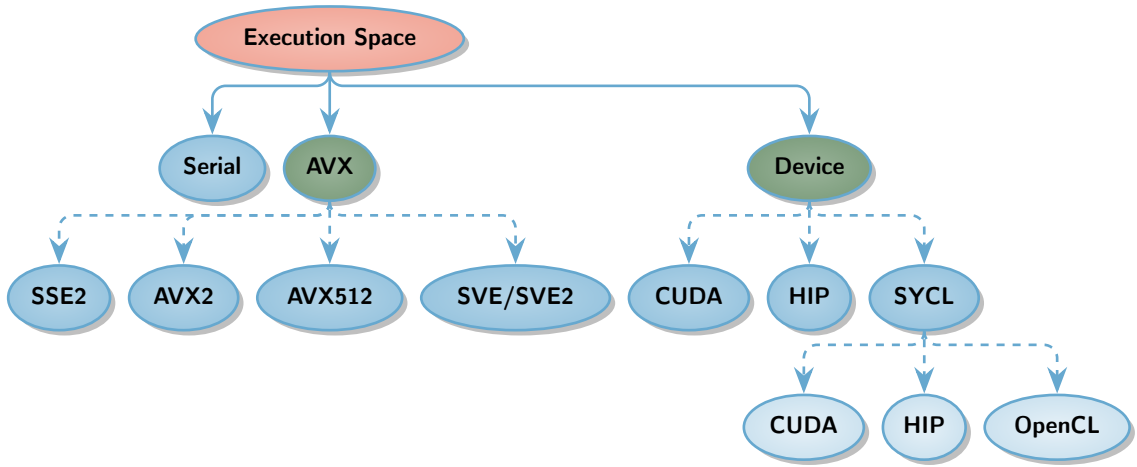

 %To support a variety of execution back-ends, notably computers with SIMD vectorisation and modern heterogeneous architectures exploiting discrete GPUs, an {\it execution space} model was introduced.
 An execution space is an intermediate representation for a family of back-end models that exploit a common parallelisation strategy. In addition, each execution space is associated with a memory space in which the algorithm is executed.  To illustrate this, the execution spaces used in Nektar++ are schematically outlined in Figure \ref{fig:execution_space} and include:
\begin{itemize}
\item \textbf{Serial}: Represents the default execution space in Nektar++ and is associated with the Host memory space.
\item \textbf{AVX}: Represents back-ends for explicit SIMD vectorisation, currently including support for SSE2, AVX2, AVX512, and SVE/SVE2. This execution space is also associated with a Host memory space.
\item \textbf{Device}: Represents back-ends for offloading to GPU accelerators, currently including support for CUDA, HIP, and SYCL. This execution space is associated with a Device memory space.
\end{itemize}

The AVX back-end maximises the potential of SIMD vectorisation, since in previous investigations \cite{moxey-2019b} we observed that changing the data storage layout to cluster data entries from four to eight elements into a vector register improves cache reuse, which is crucial to boost the performance of AVX. However, similar to the xsimd \cite{xsimd} project, the various SIMD instruction sets are unified by defining wrappers for scalar and vector data types and by overloading arithmetic operators to use different SIMD intrinsics, which allows a unified formulation for the Serial and AVX back-ends.

Similarly, the CUDA, HIP, and SYCL back-ends are unified into a single formulation by defining wrapper functions and types. The proposed strategy reduces the actual implementation into two coding descriptions, one for the Serial/AVX back-ends and one for the Device back-ends, which cover the eight back-ends shown in Figure \ref{fig:execution_space}. Within these, we also consider the different implementation strategies discussed in the next section.  

\subsection{Implementation strategies}
\label{sec:implementation:strategies}
The performance of code for a specific computing task can be improved primarily in two ways. The first is through better utilisation of the available hardware, which is managed through the selection of {\it execution space}. The second is through changes to the underlying algorithm of the task. For example, a spectral/$hp$ element Helmholtz operator can be realised into multiple equivalent, but algorithmically different, formulations as already highlighted in Section \ref{sec:detail}. In a matrix-based implementation, we pre-assemble and store the dense block matrices for each element. Alternatively, in a matrix-free implementation, we do not assemble the elemental matrices, but rather compute the entries on-the-fly during operator evaluation. Similarly, we can exploit the operator structure induced by the tensor-product basis construction in two and three dimensions; rather than perform a single matrix operation on an element, we can instead perform a series of sub-operations along 1D coordinates. Different formulations exhibit different memory access patterns and operation costs. Furthermore, there are also opportunities for loop-level optimisations and kernel fusion.

Meshing complex geometries necessitates the use of unstructured meshes with various element types. The properties of different workloads also affect the performance and should be considered in the implementation, including the shape of the elemental expansion (i.e. triangle, quadrilateral in 2D), the geometric deformation of the element (i.e. deformed or a regular affine transformation), the polynomial degree of the expansion and the integration/quadrature order, which may also vary in each elemental coordinate direction. As previously advocated in earlier work \cite{moxey-2016}, providing multiple implementations for the same operation allows improved performance to be achieved across different architectures. This study also highlights the performance benefits of treating groups of identically structured elements together where they have the same shape, polynomial degree, quadrature order, coordinate mapping type (i.e. deformed or affine), in terms of memory organisation, data access patterns and computational efficiency; in this work, we use the term \emph{block} to describe such a grouping of elemental data, in reference to its organisation, with constant stride, in memory. Different implementation choices may be made per block for a given operation.

In previous work (\cite{Cantwell2011} and \cite{nektar2020}) we have explored these aspects and identified several beneficial \emph{implementation strategies} which we discuss in the remainder of this section, and illustrate through the backward transform operator which maps elemental coefficients to physical-space quadrature points within an element, or equivalently interpolates a nodal expansion to a set of quadrature points that are not at the node points of the expansion. The operation is commonly denoted as $\boldsymbol{u}={\bm B} \boldsymbol{\hat{u}}$, and is described in Section \ref{sec:detail}. We now present four different, but equivalent in action, backward transform operation implementations (Algorithms-\ref{alg:bwdtrans-stdmat}--\ref{alg:bwdtrans-sumfactop}) and their performance characteristics, and then subsequently present two different, but equivalent in action, Helmholtz operators as outlined in Section~\ref{sec:reinterp} and their performance characteristics.

\subsubsection{Standard matrix (StdMat): with and without element grouping}
\label{sec:stdmat}
This approach performs the target operation by evaluating a matrix operation with reference-element matrices on multiple elements together within a block. This approach is possible for all execution spaces. Elements of the same block are first mapped back to the reference element where the operator is identical for all elements. Element-wise vector operations are used to apply any element-specific geometric scaling, followed by a single matrix-matrix multiplication to evaluate the operator on the reference element. The second operand consists of columns, each corresponding to degrees of freedom of a single element. This approach is matrix-free in that it does not involve creating a global matrix system over many elements. It is commonly adopted in nodal spectral/$hp$ element implementations, especially for triangles and tetrahedra, where further factorisation is impossible to achieve if the nodal expansions does not utilise a tensor product structure.
    
For the StdMat implementation of the backward transform, we perform a general matrix-matrix multiplication (GEMM) as shown in Algorithm~\ref{alg:bwdtrans-stdmat}. No element-wise scaling is required for this operator. $N_P$ is the total number of degrees of freedom in an element and $N_Q$ is the total number of quadrature points in an element. $N_e$ is the number of elements in a block. The GEMM can be performed by using a suitable BLAS implementations, such as the LIBXSMM \cite{heinecke2016libxsmm}, cuBLAS, hipBLAS, and oneMath \cite{oneMath} libraries for the AVX, CUDA, HIP, and SYCL back-ends, respectively.

\begin{algorithm}
	\caption{\emph{Standard matrix} implementation of the backward transform.}
	\label{alg:bwdtrans-stdmat}
	\begin{algorithmic}[1]
		\Require $\hat{\boldsymbol{u}}[N_P \times N_e]$; $ \bm B[N_Q \times N_{P}]$.
		\Ensure $ \bm u[N_Q \times N_e] $
        \medskip
		\State GEMM: $\bm u \gets \bm B\hat{\boldsymbol{u}}$
	\end{algorithmic}
\end{algorithm}

For GPUs, larger values of $N_e$ are preferred in order to hide memory latency and saturate the processor, so we process all elements in one GEMM call. However, for CPUs, we only process a small group of elements in each GEMM call when fusing other kernels in order to restrict the intermediate data size and improve data locality. The LIBXSMM library is optimised for small matrix-matrix multiplication ($M N K \leq 64^3$ where $M$, $K$ and $N$ are the dimensions of the matrices) and good performance is expected at lower polynomial degree. Therefore, in this approach, we operate on small sets of $W_v$ elements, as illustrated in Algorithm~\ref{alg:bwdtrans-stdmat-simd}. We add padding where necessary, when $N_e$ is not divisible by $W_v$. We found that using SIMD width of the processor as the group size $W_V$ is often a good choice to achieve better performance. To ensure data locality and efficient memory access, the data are \textit{interleaved} so that consecutive lanes within a vector group (i.e. a SIMD register) are accessing contiguous data during memory load and store. These extra data movements take time and compromise performance. To minimise the impact of interleaving, data are left in an interleaved state to avoid unnecessary deinterleaving and re-interleaving between multiple vectorised operations. 
This can also be considered as a \textit{vectorised} \emph{Standard matrix} implementation.

\begin{algorithm}
	\caption{\emph{Standard matrix} implementation of the backwards-transform operator with element grouping.}
	\label{alg:bwdtrans-stdmat-simd}
	\begin{algorithmic}[1]
		\Require $\hat{\boldsymbol{u}}[N_P \times N_e]$; $ \bm B[N_Q \times N_{P}]$.
	    \Ensure  $ \bm u[N_Q\times N_e]$
        \medskip
		\For{group id $g$ from 1 to $N_e/W_V$ }
            \State Interleave  $\hat{\bm{u} }$ of $W_V$ elements in the group g into ${\hat{\boldsymbol{u}}_g[W_V \times N_P]}$ (if required)
            %for $g=1,\dots, N_e/W_V$
            \State GEMM: $\bm u_g \gets \hat{\boldsymbol{u}}_g  \bm B^T$
        \State De-interleave  ${\bm u_g }$ back to original ordering (if required)
		\EndFor
	\end{algorithmic}
\end{algorithm}

\subsubsection{Sum-factorisation (SumFac)}
\label{sec:sumfac}

If the elemental basis is constructed as a tensor product of one-dimensional bases, one can exploit sum-factorisation to reduce the operation cost from $O(P^{2d})$ to $O(P^{d+1})$ and also reduce the matrix size from $O(P^{2d})$ to $O(P^{2})$ \cite{KarniadakisSherwin2005}, where $d$ is the space-dimension of the problem. This is particularly beneficial for higher-order spectral/$hp$ elements. Sum-factorisation decomposes the elemental operation into a series of 1D operations, performing the target operation in each coordinate direction separately. Therefore we have that $N_P=P_1 \times P_2$ and $N_Q=Q_1 \times Q_2$, where $P_d$ and $Q_d$ are the number of polynomial coefficients and quadrature values, respectively, in direction $d$. As with the standard-matrix implementation, this approach is also possible for all execution spaces.

\begin{figure}[htp]
    \centering
    \includegraphics[width=0.75\linewidth]{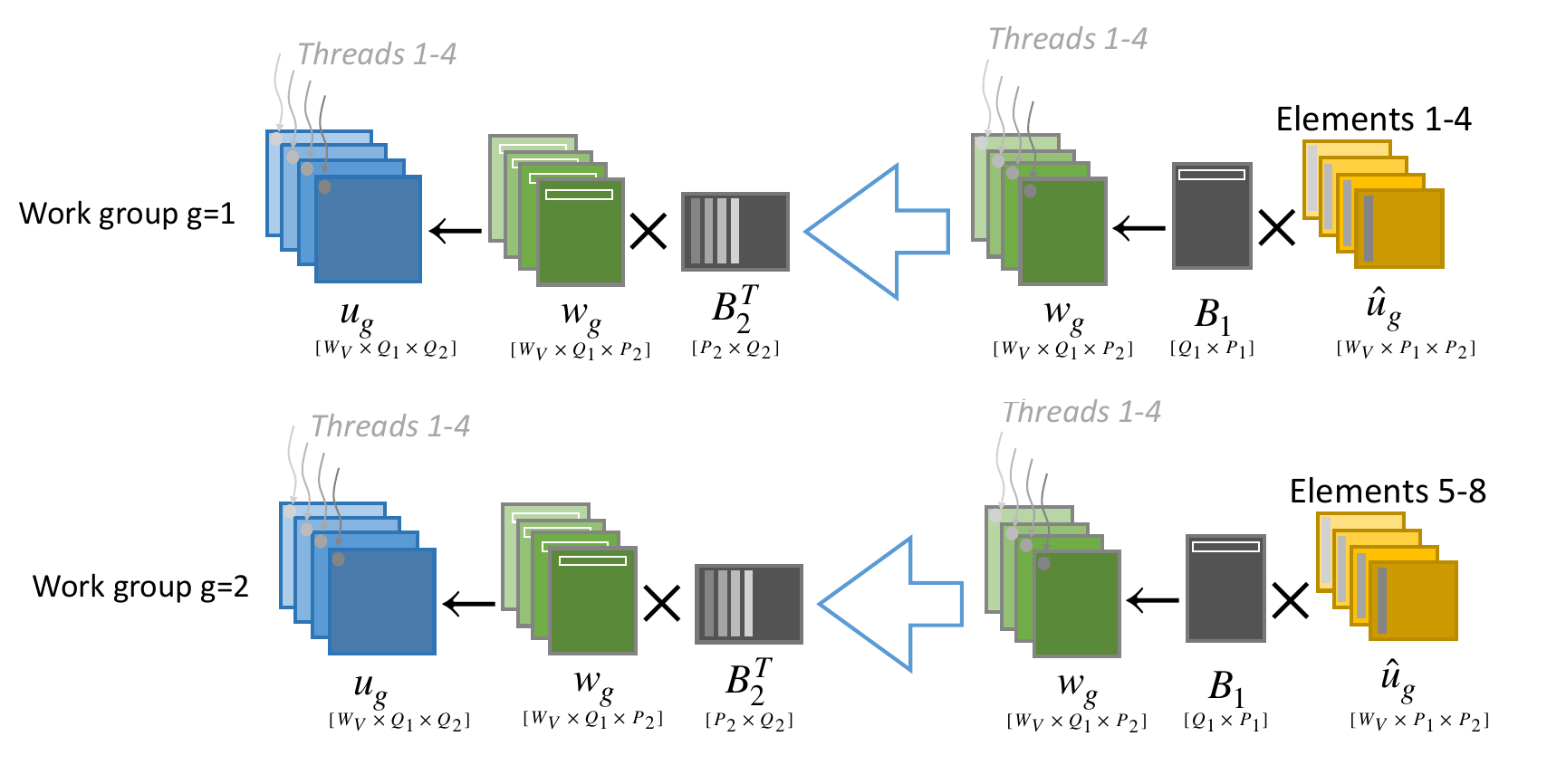}
    \caption{Diagram showing the SumFac implementation on GPU. The four grey dots are stored contiguously in the memory and accessed by different threads. For illustrative purposes, we suppose there are 4 threads in each work group. Each work group processes a different element group.}
    \label{fig:sumfac}
\end{figure}

The vectorised sum-factorisation approach is a variant of the sum-factorisation algorithm, supporting both Serial/AVX and Device back-ends. The Device back-ends introduce an additional level of parallelisation over multiple vector groups (i.e. via CUDA/HIP Grid). The pseudo-code of the sum-factorisation algorithm on quadrilateral elements is given in Algorithm~\ref{alg:bwdtrans-sumfac} and its implementation on GPU is illustrated in Figure \ref{fig:sumfac}. It is important to note that line 3 and 4  perform GEMM-like tasks, but they do not follow the standard GEMM interfaces. For elements involving simplices, such as triangles or tetrahedrons, line 4 becomes a triangular matrix multiplication rather than GEMM. In Nektar++, these operations are manually coded rather than utilising BLAS routines, which also allows us to fuse additional pre and post operations if required.

\begin{algorithm}
	\caption{\emph{Sum-factorisation} (SumFac) algorithm for the backwards-transform on quadrilateral elements.}
	\label{alg:bwdtrans-sumfac}
	\begin{algorithmic}[1]
		\Require $\hat{\boldsymbol{u}}[N_P \times N_e]$; $\bm B_d[Q_d\times P_d] , d=1,2$. 
		\Ensure $ \bm u[N_Q\times N_e]$.        
        \medskip
		\For{group id $g$ from 1 to $N_e/W_V$ } \Comment{Parallel processing  multiple groups if possible}
		       \State Interleave  $\hat{\bm{u} }$ of $W_V$ elements in the group g into ${\hat{\boldsymbol{u}}_g[W_V \times P_1 \times P_2]}$ (if required)
		\State GEMM and transpose: $\bm w_g[W_V \times Q_1 \times P_2] \gets \bm B_1 \hat{\boldsymbol{u}}_g $ \Comment{Parallel processing $W_V$ entries}
		\State GEMM and transpose: $\bm u_g[W_V \times Q_1\times Q_2] \gets  \boldsymbol{w}_g  \bm B_2^T $ \Comment{Parallel processing $W_V$ entries}
	  \State De-interleave  ${\bm u_g }$ back  to original ordering (if required)
		\EndFor
	\end{algorithmic}
\end{algorithm}

%For Serial/AVX back-ends, parallelisation only happens on the $W_V$ level, falling back to Serial execution when $W_V=1$.

\subsubsection{Sum-factorisation: Threaded on Output-Point (SumFacTOP)}
\label{sec:sumfactop}
This is a variant of the sum-factorisation implementation designed for threading type operations that are used in the Device execution space. In the sum-factorisation approach of Section~\ref{sec:sumfac}, each element is processed by a single thread on the device. In SumFacTOP, each element is processed by one work group (i.e. CUDA/HIP thread block) to exploit the multi-level hierarchical parallelism of GPUs as previously outlined in \cite{Eichstadt2020}. Each thread within a work group processes one output degree-of-freedom, resulting in an algorithm that exploits parallelism both across and within each spectral/$hp$ element.

\begin{algorithm}
	\caption{\emph{Sum-factorisation threaded on output point} (SumFacTOP) implementation of the backwards-transform on quadrilateral elements.}
	\label{alg:bwdtrans-sumfactop}
	\begin{algorithmic}[1]
		\Require $\hat{\boldsymbol{u}}[N_P \times N_e]$; $\bm B_d[Q_d\times P_d], d=1,2$. 
		\Ensure $ \bm u[N_Q\times N_e]$.
		% \State \textcolor{blue}{// Loop over vector group ()}
        \medskip
		\For{element id $e$ from 1 to $N_e$ } \Comment{Parallel processing multiple elements via multiple GPU work groups}
        \State Retrieve input of element e $ \hat{\bm u}_e [P_1\times P_2]$ from $\hat{\bm u}$
        \For{value $idx$ from 1 to $Q_1 \times P_2$ } \Comment{Parallel process entries in $\bm w_e$ over GPU threads}
    	    \State $i = idx/P_2$, $q = idx\,\%\,P_2$
            \State DOT: $\bm w_e[i\times P_2 + q] \gets \bm B_1[i][:]  \cdot\hat{\boldsymbol{u}}_e[:] [q]$ 
        \EndFor
        \For{value $idx$ from 1 to  $Q_1 \times Q_2$ } \Comment{Parallel process entries in $\bm u_e$ over GPU threads}
    	    \State $j = idx\,\%\,Q_2$, $i = idx/Q_2$
            \State DOT: $\bm u_e[i\times Q_2 + j] \gets \bm w_e[i][:] \cdot \bm B^T_2[:] [j]$ 
        \EndFor
        \State Write output of element e $\bm u_e [Q_1\times Q_2]$ back to $\bm u$
		%\State Epilogue operations on ${\bm u_e }$ (if any)
        %\Comment{Parallel processing $Q_1\times Q_2$ entries in $\boldsymbol{u}_e$}
		\EndFor
	\end{algorithmic}
\end{algorithm}

\begin{figure}[htb]
    \centering
    \includegraphics[width=0.75\linewidth]{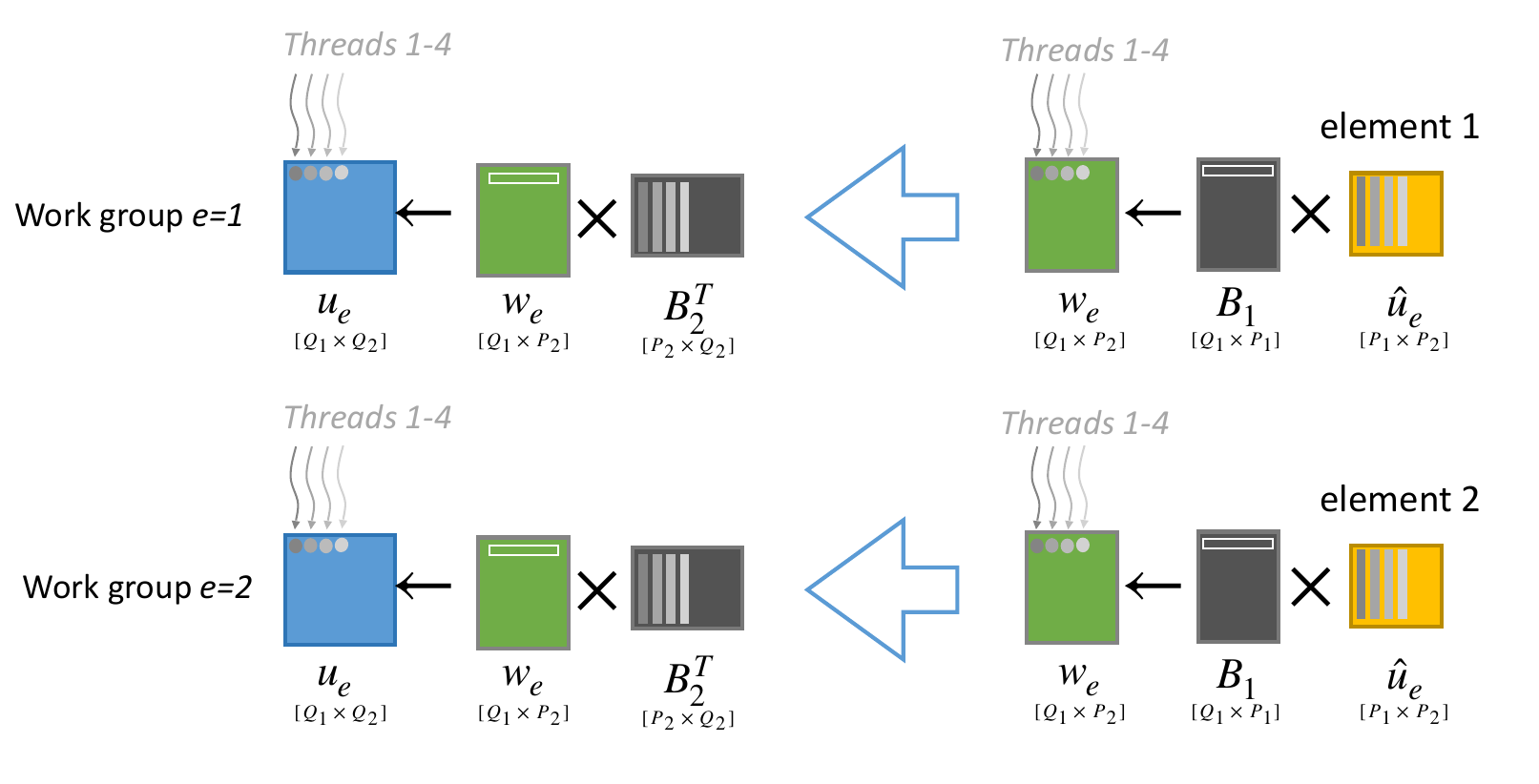}
    \caption{Diagram showing the SumFacTOP implementation on GPU. The four grey dots are stored contiguously in memory and accessed by different threads. For illustrative purposes, we suppose there are four threads in each work group. Each work group processes a different element.}
    \label{fig:sumfactop}
\end{figure}

The pseudo-code of the SumFacTOP algorithm is given in Algorithm~\ref{alg:bwdtrans-sumfactop} and its implementation on GPU is illustrated in Figure \ref{fig:sumfactop}. In this algorithm we use the percent symbol to denote modulus operator. To ensure efficient memory access, the appropriate part of the input vector, $\hat{\boldsymbol{u}}_e[P_1 \times P_2 ]$, and the basis data, $\bm B_1$ and $\bm B_2$, should be first loaded into shared memory. This should also be used to store intermediate results, such as $\bm w_e[Q_1 \times P_2]$. The innermost level of parallelism exposes sufficient work when elements have a medium to high polynomial degree. For instance, a quadrilateral element of degree six, which has $7^2$ DOFs and $8^2$ quadrarture points, or a hexahedral element of degree two, which has $3^3$ DOFs and $4^3$ quadrature points, has enough output quadrature points to fill several CUDA/HIP warps/wavefronts, each of which contains 32/64 threads, respectively.

%%
%A schematic representation of the implementation model is shown in Figure \ref{fig:implementation_model}. More detailed explanations and examples of each implementation model are also presented in Section \ref{subsec:algo}. 
%%
%\begin{figure}[htb]
%\centering
%\begin{forest}
%  for tree={
%      font=\sffamily\bfseries,
%      line width=1pt,
%      draw=linecol,
%      ellip,
%      align=center,
%      child anchor=north,
%      parent anchor=south,
%      drop shadow,
%      l sep+=12.5pt,
%      edge path={
%        \noexpand\path[color=linecol, rounded corners=5pt, >={Stealth[length=10pt]}, line width=1pt, ->, \forestoption{edge}]
%          (!u.parent anchor) -- +(0,-5pt) -|
%          (.child anchor)\forestoption{edge label};
%        },
%  }
%  [Implementation, inner color=col1in, outer color=col1out, calign=child, calign child=2
%    [StdMat, inner color=col2in, outer color=col2out, edge=dashed]
%    [SumFac, inner color=col2in, outer color=col2out, edge=dashed]
%    [SumFacTOP, inner color=col2in, outer color=col2out, edge=dashed]
%  ]
%\end{forest}
%\caption{Nektar++ implementation models.}
%\label{fig:implementation_model}
%\end{figure}
%%%

\subsubsection{Alternative implementations of the Helmholtz operator} \label{sec:helm-impl}

Section~\ref{sec:reinterp} outlined two mathematically equivalent formulations for the Helmholtz operator. These lead to different implementations in the code and consequently result in  different performance. From Equation (\ref{eq:helm}) and Equation (\ref{eq:helm1}), we can derive two distinct designs as shown in Algorithm~\ref{alg:helm_design1} and Algorithm~\ref{alg:helm_design2}.
% %
In these algorithms, we have used the following shorthand: BwdTrans to denote a backwards transform from the coefficient space of the polynomial expansion to the quadrature points; StdDeriv denotes the derivative in the standard element; SumStdDeriv is a StdDeriv using sum-factorisation; IProductWRTBase denotes an inner product with respect to the expansion space and IProductWRTDerivBase denote the inner product with respect to the derivative of the expansion space. 

\begin{algorithm}
\caption{Helmholtz implementation by non-collocated approach}
\label{alg:helm_design1}
\begin{algorithmic}[1]
\Require $\hat{\boldsymbol{u}}$, $\bm{D}_{\xi_d}\bm{B}$, $\lambda$, $\frac{\partial \xi_i}{\partial x_j}$ and $\omega J$ at quadrature points.
\Ensure  LHS results.
    \State BwdTrans:           $\boldsymbol{u} \gets \bm{B}\hat{\boldsymbol{u}}$
    \State BwdTrans \& StdDeriv:   $\boldsymbol{v}_1 \gets (\bm{D}_{\xi_1}\bm{B}) \hat{\boldsymbol{u}},\, \boldsymbol{v}_2 \gets (\bm{D}_{\xi_2} \bm{B}) \hat{\boldsymbol{u}}$
    \State ApplyMetric: $ \boldsymbol{v}_1' \gets \bm{\Lambda}_{11}\boldsymbol{v}_1 + \bm{\Lambda}_{12}\boldsymbol{v}_2,\, \boldsymbol{v}_2' \gets \bm{\Lambda}_{21}\boldsymbol{v}_1 + \bm{\Lambda}_{22}\boldsymbol{v}_2$, 
        \Comment{$\bm{\Lambda}$ consists of $\frac{\partial \xi_i}{\partial x_j}$. See Equation (\ref{eq:helm1}).}
    \State IProductWRTDerivBase:   $\text{LHS} \gets (\bm{D}_{\xi_1}\bm{B})^{T} \boldsymbol{W} \boldsymbol{v}'_1 + (\bm{D}_{\xi_2}\bm{B})^{T} \boldsymbol{W}\boldsymbol{v}'_2 $
    \State IProductWRTBase: $\text{LHS} \gets \text{LHS} + \lambda\bm{B}^T \boldsymbol{W} \boldsymbol{u}$
\end{algorithmic}
\end{algorithm}

\begin{algorithm}
\caption{Helmholtz implementation by collocated approach}
\label{alg:helm_design2}
\begin{algorithmic}[1]
\Require $\hat{\boldsymbol{u}}$, $\bm{B}$, $\bm{D}_{\xi_d}$, $\lambda$, $\frac{\partial \xi_i}{\partial x_j}$ and $\omega J$ at quadrature points.
\Ensure  LHS results.
    \State BwdTrans: $\boldsymbol{u} \gets \bm{B}\hat{\boldsymbol{u}}$
    \State StdDeriv:      $\boldsymbol{v}_1 \gets \bm{D}_{\xi_1} \boldsymbol{u},\, \boldsymbol{v}_2 \gets \bm{D}_{\xi_2} \boldsymbol{u}$
    \State ApplyMetric: $ \boldsymbol{v}_1' \gets \boldsymbol{W}(\bm{\Lambda}_{11}\boldsymbol{v}_1 + \bm{\Lambda}_{12}\boldsymbol{v}_2),\, \boldsymbol{v}_2' \gets \boldsymbol{W}(\bm{\Lambda}_{21}\boldsymbol{v}_1 + \bm{\Lambda}_{22}\boldsymbol{v}_2)$, 
        \Comment{$\bm{\Lambda}$ consists of $\frac{\partial \xi_i}{\partial x_j}$. See Equation (\ref{eq:helm1}).}
    \State SumStdDeriv:   $\boldsymbol{u}'\gets \bm{D}^T_{\xi_1} \boldsymbol{v}'_1 + \bm{D}^T_{\xi_2} \boldsymbol{v}'_2 + \lambda \boldsymbol{W}\boldsymbol{u}$
    \State Restrict: $\text{LHS} \gets \bm{B}^T \boldsymbol{u}'$
\end{algorithmic}
\end{algorithm}

We observe that for a quadrilateral element $\bm{B}$ and $\bm{D}_{\xi_d} \bm{B}$ are dense matrices of the same shape $Q^2 \times P^2$, but can be optimised by sum-factorisation, reducing the matrix-vector multiplication cost from $Q^2P^2$ to $P^2Q + PQ^2$. The reduction is even more substantial in three dimensions. $\bm{D}_{\xi_d}$ in its most general form is a square matrix of size $Q^2 \times Q^2$, but with the collocation and sum-factorisation properties of tensor Lagrange product bases, it can be regarded as a blocked matrix with a matrix-vector multiplication cost of $Q^3$. In the SumFac/SumFacTOP implementation, the non-collocated approach requires six sum-factorisation operations, while the collocated approach requires two sum-factorisation operations and four blocked-matrix operations, which may have lower floating-point operation cost. Our benchmark results given in the next section confirm that the collocated approach indeed notably outperforms the non-collocated approach in the case $Q=P+1$ for all the tested orders and element types. However, for the StdMat implementation, the  $\bm{D}_{\xi_d}\bm B$ matrix is dense but the dimension is typically smaller than the $\bm{D}_{\xi_d}$ matrix and so  the non-collocated approach almost always has lower floating point cost and data loads, and is therefore adopted  in the StdMat implementation.

We next explore the performance of the different implementation approaches outlined above across different three-dimensional element shapes and polynomial orders.

%% file: Results.tex
\section{Results}
\label{sec:results}
In benchmarking the different implementation strategies discussed in Section \ref{sec:implementation}, we have leveraged some of the benchmark test cases proposed under the CEED project \cite{brown2021libceed}. In the following, we present measurements for the BK1 and BK3 cases which consider the throughput of the mass matrix and stiffness/Helmholtz matrix operations as a function of the total number of elemental degrees of freedom on different nodes, with either GPU and AVX acceleration. Considering the action ${\bm f} = {\bm{M v}}$, these tests therefore measure how many output degrees of freedom of ${\bm f}$ can be computed per second when calculating the action of a matrix $\bm M$ (i.e. the mass or Helmholtz matrix) on an input vector $\bm v$. 

In the following tests, we consider performance on an NVIDIA GH200 Grace Hopper Superchip and an Intel Xeon 6526Y CPU with 32 cores and AVX-512 instructions. Additional results for a larger CPU (AMD EPYC 9554) can be found in Appendix~\ref{app:results:epyc}. For different polynomial degrees, we have tested the standard matrix (StdMat) approach, the sum-factorisation (SumFac) using a vectorised approach, and the sum-factorisation approach threaded on output point (SumFacTOP). We first highlight some general observations relating to these tests: 
\begin{itemize}
    \item Consistent with the paper by the CEED team \cite{fischer-2020}, all tests were performed assuming the elements were curvilinear, approximated by the same degree of polynomial as the solution with iso-parametric mappings and were not allowed to exploit the global tensor-product structure of the element layout in the mesh. Similarly to \cite{fischer-2020}, we have also precomputed the geometric matrices.
    \item In the CEED project, we are aware of results for hexahedral elements, but we have not been able to find published measurements for high-order prismatic, pyramidic and tetrahedral elements that we present below. 
    \item The complexity of different operations gives rise to different optimal performance from the three implementation strategies outlined in Section~\ref{sec:implementation:strategies}. In general, the algorithmically simpler mass matrix operation benefits from the StdMat implementation, whereas the more algorithmically complex Helmholtz matrix operation benefits from the SumFac or SumFacTOP type implementation, depending on the architecture. 
    %\item On the Intel Xeon 6526Y AVX node, we observe relatively consistent throughput performance over all shape types (i.e. hexahedral (Hex), prismatic (Prism), pyramidic (Pyr) and tetrahedral (Tet) elements). However, on the NVIDIA GH200 Grace Hopper GPU node we observe a reduction in throughput as we move from hexahedral to tetrahedral elements which might be expected as there is a greater algorithmic intensity per output degree of freedom in the tetrahedral high-order element as compared to that of a hexahedral element. However, we only observe a reduction of at most a factor of $2.5$ where as the additional intensity of the operator is around a factor of six. 
    For the latter, we observe a reduction in throughput as we move from hexahedral to tetrahedral elements, which is likely to be predominantly caused by the larger memory requirement for the basis data in tetrahedra. For example, on the NVIDIA GH200 Grace Hopper Superchip, we observe a reduction of at most a factor of $2.5$, whereas the computational cost of the operator is around a factor of six higher, owing to the tensor-product collapsed coordinate system used in tetrahedral elements. 
\end{itemize}

\subsection{NVIDIA GH200 Grace Hopper Superchip}
\begin{figure}[htbp]
\includegraphics[width=\textwidth]{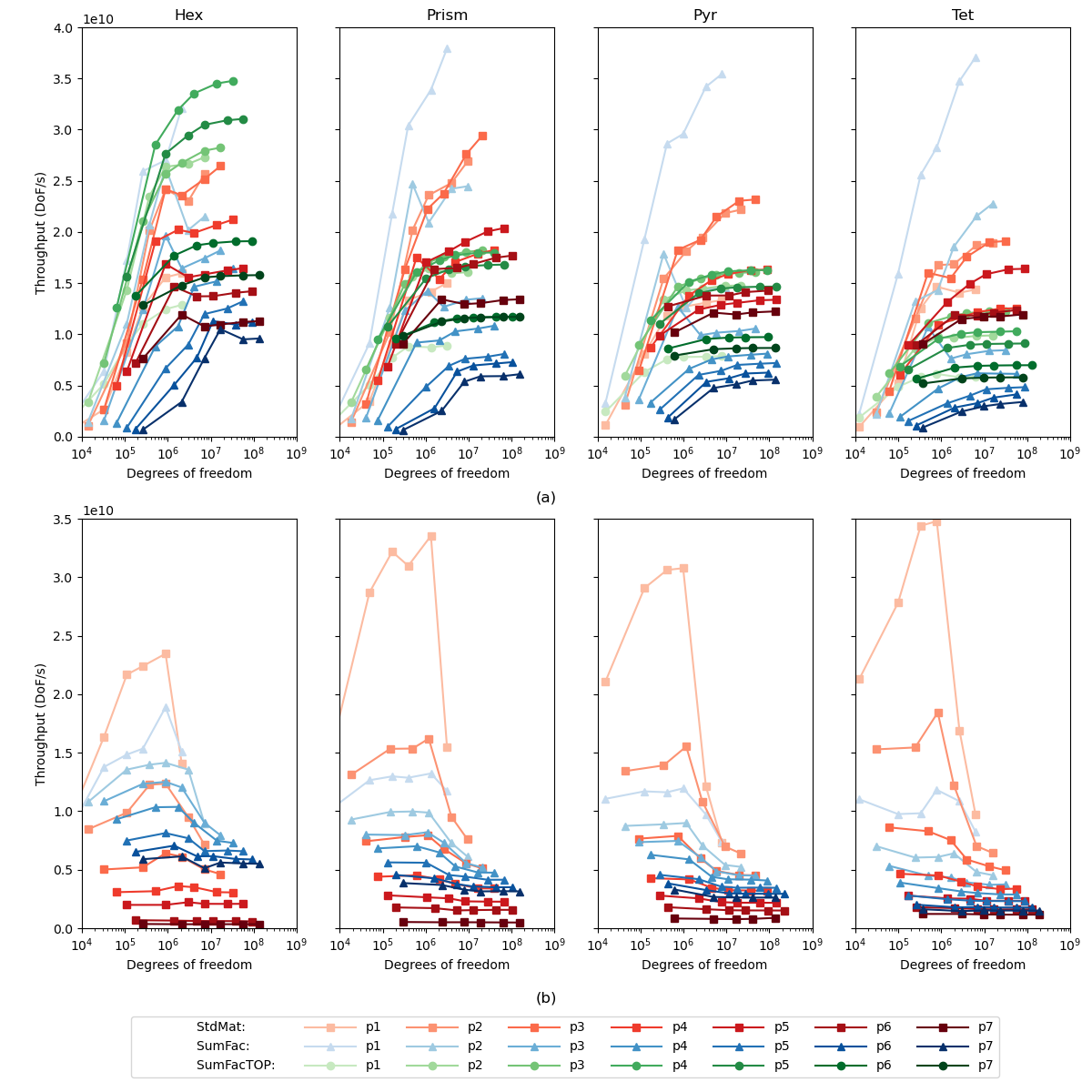}
\caption{(a) NVIDIA GH200 Grace Hopper Superchip and (b) Intel Xeon 6526Y throughput performance versus elemental degrees of freedom for the mass operator for hexahedral (Hex), prismatic (Prism), pyramidic (Pyr) and tetrahedral (Tet) elements. For different polynomial degree (P) the standard matrix (StdMat) approach is labelled in red, the vectorised sum-factorisation (SumFac) is labelled in blue and the sum-factorisation threaded on output point (SumFacTOP) is labelled in green. }
\label{fig:results_mass}
\end{figure}

\begin{figure}[htbp]
\includegraphics[width=\textwidth]{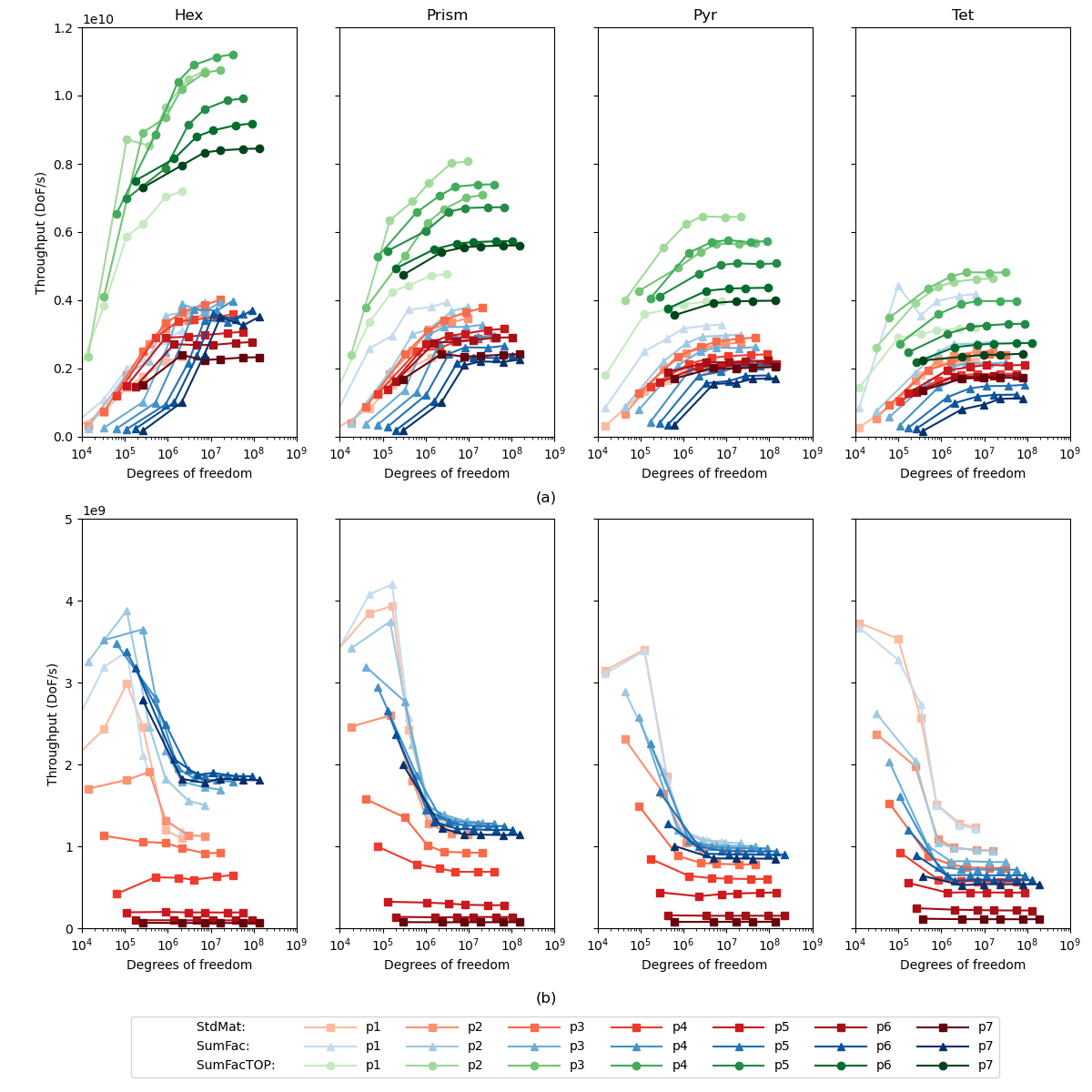}
\caption{(a) NVIDIA GH200 Grace Hopper Superchip  and (b) Intel Xeon 6526Y  throughput performance versus elemental degrees of freedom for the Helmholtz operator for hexahedral (Hex), prismatic (Prism), pyramidic (Pyr), and tetrahedral (Tet) elements. For different polynomial degree (P) the standard matrix (StdMat) approach is labelled in red, the vectorised sum-factorisation (SumFac) is labelled in blue and the sum-factorisation threaded on output point (SumFacTOP) is labelled in green. }
\label{fig:results_helm}
\end{figure}

The results obtained with the CUDA back-end on a NVIDIA GH200 Grace Hopper Superchip are shown in Figures~\ref{fig:results_mass}(a) and \ref{fig:results_helm}(a). In these plots we show the throughput performance for the mass and Helmholtz matrix operations, as a function of the elemental degrees of freedom, on hexahedral (Hex), prismatic (Prism), pyramidic (Pyr), and tetrahedral (Tet) elements for a range of polynomial degrees (P). In these plots we compare three implementation strategies namely: the standard matrix (StdMat) approach (denoted in red), outlined in section \ref{sec:stdmat}; the vectorised sum-factorisation (SumFac) approach (denoted in blue), outlined in section \ref{sec:sumfac}, and; the sum-factorisation threaded-on-output-point (SumFacTOP) (denoted in green), outlined in section \ref{sec:sumfactop}.

First considering the less algorithmically intensive mass operator in Figure~\ref{fig:results_mass}(a) we first observe that the StdMat (red lines) and SumFacTOP (green lines) tend to provide improved throughput performance for the higher polynomial degree above $P=2$. However, the vectorised sum-factorisation (SumFac) implementation, highlighted by blue lines, shows excellent  performance for shapes involving simplexes faces (i.e. prisms, pyramids, and tetrahedra) for low polynomial degree with $P = 1$ and $P = 2$. This is likely due to the small number of degrees of freedom per element, which allows efficient use of registers. Nevertheless, the performance of this SumFac implementation quickly drops as the polynomial degree increases, with a markedly lower performance beyond $P = 4$. However, the SumFacTOP implementation shows better performance compared to the StdMat implementation for hexahedral elements, while the opposite is observed for simplex elements. Finally, SumFacTOP typically shows low performance for $P = 1$ as for such a low-degree polynomial, there are insufficient degrees of freedom per element to fully utilise a CUDA warp of 32 lanes. 

For the high-arithmetic intensity Helmholtz operator in Figure~\ref{fig:results_helm}(a), the SumFacTOP implementation outperforms the StdMat implementation across all polynomial degrees, shapes, and degrees of freedom. This implementation notably allows the use of shared memory for intermediate calculations, allowing efficient kernel fusing and data re-use on an element basis by exploiting the hierarchical parallelism characteristics of GPUs. On the other hand, both the SumFac and StdMat implementations involve multiple load/store from global memory. 

The results also indicate that elements with simplex faces (prisms, pyramids, and tetrahedra) offer lower performance compared to hexahedral elements for both operators due to the cost of integrating using one-dimensional quadrature based on the Duffy transformation (see Figure~\ref{fig:collapsed}). This tensor product integration approach notably introduces additional degrees of freedom in the quadrature space compared to the number of expansion degrees of freedom. Furthermore, building the expansion basis requires larger basis data for the warped tensor product in these elements \cite{KarniadakisSherwin2005}.  Nevertheless, whilst there is a six-fold increase in the number of quadrature points relative to the number of expansion coefficients in tetrahedral elements compared to the hexahedral elements, there is only, at most, a reduction by a factor of 2.5 in the throughput for the Helmholtz operator in Figure~\ref{fig:results_helm}(a). 

Finally, to compare against the performance achieved in the CEED project \cite{brown2021libceed} we have compiled the MFEM code on the same hardware and performed the same operation on hexahedral elements. For the mass matrix operator utilising the SumFacTOP approach on a GH200 we observe a peak throughput of $2.5\times 10^{10}$ with MFEM. This is similar to our SumFacTOP performance shown in Figure~\ref{fig:results_mass}(a) although we do observe a maximum performance of $3.5\times 10^{10}$.  Similarly for the Helmholtz operator, the MFEM stiffness operator achieves a peak performance of $1.1 \times 10^{10}$ throughput, which is very similar to our peak performance shown in Figure~\ref{fig:results_helm}(a).  

\subsection{Intel Xeon 6526Y CPU}
Figures~\ref{fig:results_mass}(b) and \ref{fig:results_helm}(b) show the throughput performance we obtained on an Intel Xeon 6526Y dual-socket CPU node with a total of 32 cores. Once again in these plots we show the throughput performance for the mass and Helmholtz matrix operations as a function of the elemental degrees of freedom on hexahedral (Hex), prismatic (Prism), pyramidic (Pyr) and tetrahedral (Tet) elements, and for a range of polynomial degrees (P). In these plots we compare only the two relevant implementation strategies namely: the standard matrix (StdMat) approach and the vectorised sum-factorisation (SumFac) approach. The  sum-factorisation threaded-on-output-point (SumFacTop) algorithm is not possible on this architecture. 

Considering first the lower algorithmic intensity mass matrix operator in Figure~\ref{fig:results_mass}(b) we observe that at lower polynomial degrees the StdMat (red lines) implementation performs better on all elemental shapes. However, as we increase the polynomial degree, the performance of this implementation drops off dramatically and the vectorised SumFac (blue lines) algorithm performs better for $P>3$. This drop-off of performance is expected, both from the fact that the LIBXSMM library is specifically optimized from small matrix multiplication, as well as the better algorithmic efficiency of the SumFac approach at higher polynomial degree. However, the small number of degrees of freedom in the tetrahedral expansion leads to the StdMat approach remaining competitive for all polynomial degrees investigated.   For higher number of degrees of freedom at a fixed polynomial degree, which corresponds to increasing the number of elements, there is a peak in throughput performance followed by a reduction as we continue to increase the number of elements/degrees of freedom. This can be attributed to saturation of the processor cache. %Finally it is interesting to observe that the performance over the different element types is relatively consistent for all element shapes. 

Next, considering the higher algorithmic intensity Helmholtz matrix operators shown in Figure~\ref{fig:results_helm}(b), the vectorised sum-factorisation (SumFac) approach (blue lines) is markedly superior to the standard matrix (StdMat) implementation (red lines) for almost all polynomial degrees and elemental shapes. For this operator we observe a much more noticeable  performance drop-off  at high numbers of degrees of freedom which, as with the mass matrix, is likely explained by cache saturation. However, in this case there is a substantial increase in memory movement involving six geometric factors as well as the Jacobian, the later is also required for the mass matrix operator. Finally, we once again observe relatively consistent performance over all the element shapes, unlike the performance in the GH200 GPU tests.  

\subsection{Affine versus deformed elements}
Although the CEED test cases we have so far presented require that the discretisation be equivalent to a fully deformed curvilinear geometry, in our investigation it was also possible to consider the case where the elements can be geometrically described by an affine transformation. An affine transformation in a hexahedral element would mean that it maintains a parallelepiped shape where opposite faces are parallel. In this case the Jacobian and geometrical factors are independent of spatial coordinate allowing for fewer loads of geometrical information when applying the operators. For the case of a straight-edged tetrahedron, the geometry always has an affine mapping, which implies that such an optimization can have a significant effect in a general unstructured mesh. In our context, we differentiate these two cases using the terms \textit{Regular} and \textit{Deformed}, and provide slightly different implementations for them.

\begin{figure}[htb]
\includegraphics[width=\textwidth]{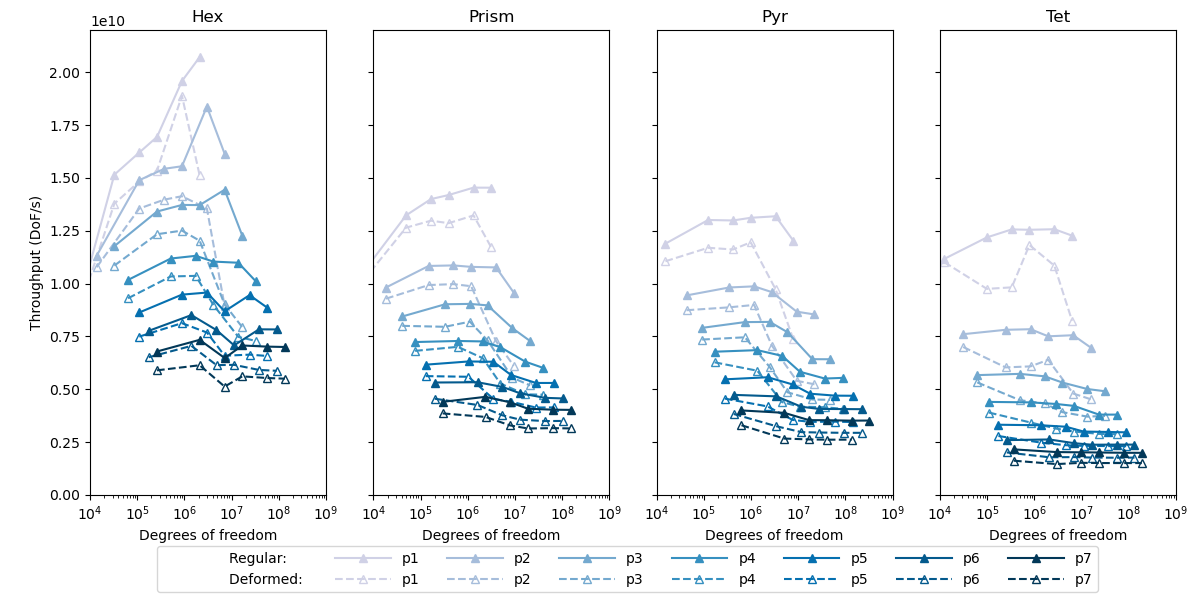}
\caption{Intel Xeon 6526Y throughput performance versus elemental degrees of freedom on mass operator for hexahedral (Hex), prismatic (Prism), pyramidic (Pyr) and tetrahedral (Tet) elements using the vectorised sum-factorisation (SumFac) implementation. The solid lines indicate throughput assuming regular (Regular) elements that have an affine transformation where as the dashed lines indicate the throughput for curvlinear (Deformed) elements. }
\label{fig:reg_def_xeon_mass}
\end{figure}

Repeating the previous tests under this affine-mapping assumption, the mass matrix operator, where there is only the Jacobian information that is now treated as a single scalar per element, exhibits a relatively consistent throughput behaviour for the Xeon nodes as shown in Figure~\ref{fig:reg_def_xeon_mass}. As the number of degrees of freedom increases, the deformed element cases the throughput drops off before that of the affine element cases, while virtually no difference is observed between the deformed and affine cases on the GH200 (see additional results in Appendix~\ref{app:results:affine}, Figure~\ref{fig:reg_def_gh200_mass}).

\begin{figure}[htb]
\includegraphics[width=\textwidth]{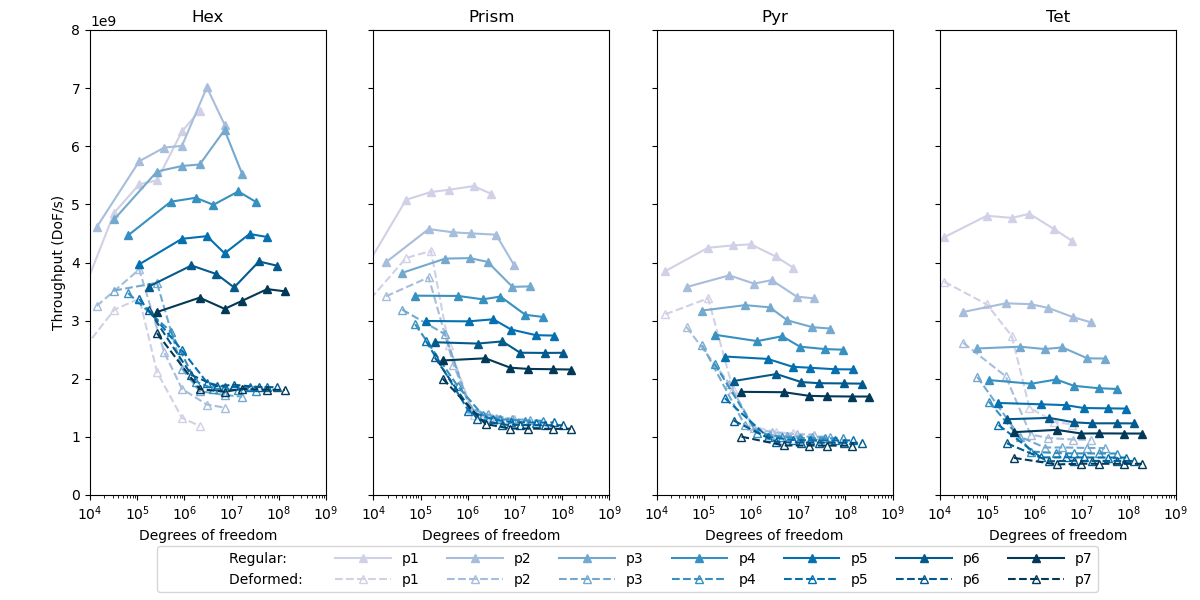}
\caption{Intel Xeon 6526Y throughput performance versus elemental degrees of freedom on Helmholtz operator for hexahedral (Hex), prismatic (Prism), pyramidic (Pyr) and tetrahedral (Tet) elements using the vectorised sum-factorisation (SumFac) implementation. The solid lines indicate throughput assuming regular (Regular) elements that have an affine transformation where as the dashed lines indicate the throughput for curvlinear (Deformed) elements. }
\label{fig:reg_def_xeon_helm}
\end{figure}

For the Helmholtz matrix operator the large size of geometrical data have a pronounced effect on whether the element is curvilinearly deformed or affine. On the GH200 GPU node, we observe quite notable throughput benefits, where the lower polynomial degrees can have up to twice the throughput (see additional results in Appendix~\ref{app:results:affine}, Figure~\ref{fig:reg_def_gh200_helm}). However, for the higher polynomial degrees tested ($P=5,6,7$) the performance was relatively similar across all elemental shapes. For the Xeon node we observe quite different performance over a range of polynomial degrees as shown in Figure~\ref{fig:reg_def_xeon_helm}. For all elemental shapes, when exploiting the affine mapping, their performance does not suffer from the notable drop-off at higher degrees of freedom, which is consistent with being able to maintain a better memory bandwidth to cache size. Indeed, the Xeon processor has an L3 cache of 1.5 MB per core whereas the GH200 node has 144GB fast HBM3e memory which is probably better able to keep the cores supplied with data.  For the Xeon node the throughput performance is roughly twice as high for the regular/affine elements as compared to the curvilinearly deformed elements. 

\subsection{Maximising collocation-based operations}
\label{sec:coll}
As detailed in section \ref{sec:reinterp} and section \ref{sec:helm-impl}, it is possible to re-interpret operations such as the inner product with respect to the derivative of the expansion basis, that arises in the weak Laplacian operator, such that we can maximally perform operations at the quadrature points. We recall that the modal expansions we adopted in all shapes can be expressed in terms of Lagrange polynomials through the quadrature points where we can maximally make use of the collocation property of the Lagrange polynomial. We then finalise the operation by restricting the evaluation back to the modal expansion. 

\begin{figure}[htb]
\includegraphics[width=\textwidth]{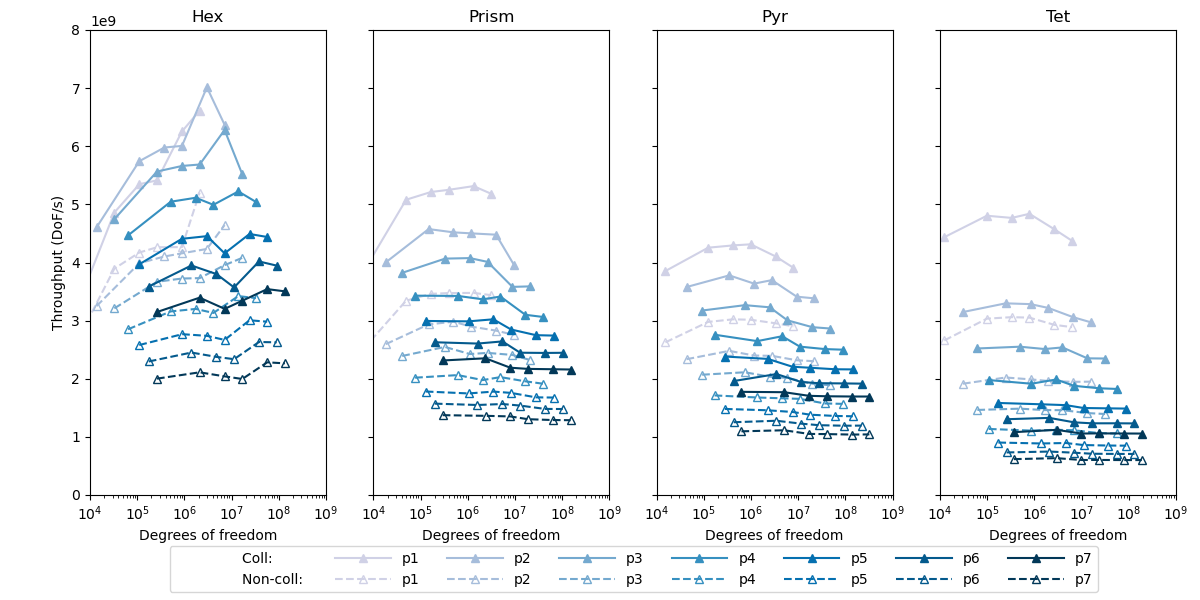}
\caption{Intel Xeon 6526Y throughput performance versus elemental degrees of freedom on Helmholtz operator for hexahedral (Hex), prismatic (Prism), pyramidic (Pyr), and tetrahedral (Tet) elements using the vectorised sum-factorisation (SumFac) implementation. The solid lines indicate throughput assuming collocated (coll) approach on regular elements where as the dashed lines indicate the throughput for the traditional non-collocated (non-coll) approach. }
\label{fig:xeon_collocated}
\end{figure}

Figure~\ref{fig:xeon_collocated} shows the results of the Helmholtz operator on all shapes on the Intel Xeon 6526Y processor. For this case, we have considered regular elements to remove the influence of the data movement associated with the geometric factors. The solid lines in this figure indicate the new approach where we are maximally using the collocated property (coll) and the dashed lines denote the previous approach where we are not maximally using the collocation property (non-coll).  For all element shapes and polynomial degrees we observe an improvement in the throughput from adopting this collocation maximisation approach. Even for the hexahedral elements, we achieve between a 25\%-50\% percent improvement as we move from $P=1$ to $P=7$. For the tetrahedral elements we observe a more consistent improvement of approximately 50\% because sum-factorisation operations on non-standard tensorial elements have lower performance and better performance is expected if they can be replaced by collocated operations on standard tensorial structures.

%% file: Summary.tex
\section{Summary}
\label{sec:summary}
Leveraging previous efforts \cite{brown2021libceed} to optimise spectral element/high order finite element methods on GPUs and vectorised CPUs,  we have 
extended this work to mixed-shaped element meshes involving not only hexahedral, but also prismatic, pyramidic, and tetrahedral shapes using tensorial expansions. Consideration of the basis in which we construct these operations highlights that the use of tensor product Lagrange base expansions can be efficiently utilized in elements of all shape types. Such an approach can  maximise operations using the collocation properties of the nodal tensorial expansion associated with classical quadrature rules.

Our study, corresponding to the BK1 and BK3 test cases of the CEED project \cite{brown2021libceed}, has  demonstrated that different operators, such as the mass and Helmholtz matrices, require alternate implementation strategies depending on the element-shape, polynomial order, and architecture type to achieve optimal performance. Three different implementation strategies were considered including: 
StdMat, utilising matrix operators in and standard element region; SumFac, utilising sum-factorisation with vectorisation, or; SumFacTOP, also utilising sum-factorisation but threaded on the output point. In general, the algorithmically simpler mass matrix operation benefits from the StdMat implementation, whereas the more algorithmically complex Helmholtz matrix operation benefits from the SumFac or SumFacTOP implementations, depending on the architecture. 

Our GPU performance tests demonstrate that the throughput of the Helmholtz operator on tetrahedral elements is at most $2.5$ times slower than on hexahedral elements, despite a ratio of six when comparing a tetrahedral element to a hexahedral element in terms of floating-point operations (flops). 
%On the Xeon AVX node, we surprisingly observe relatively consistent throughput performance over all shape types (i.e. hexahedral (Hex), prismatic (Prism), pyramidic (Pyr) and tetrahedral (Tet) elements). 

The present work has also outlined the redesign strategy proposed to extend the original Nektar++ \cite{nektar2015,nektar2020} from a purely CPU framework with some support of SIMD acceleration to a new comprehensive framework supporting both CPUs with SIMD vectorization and heterogeneous architectures with discrete GPU devices. The redesign effort notably addressed memory management, data structure, and architecture-specific optimization strategies.

\section{Acknowledgements}
We would like to acknowledge insightful conversations with Edward Erasmie-Jones, Henrik W\"ustenberg and Junjie Ye.
%and the contributions of Diego B. Renner to the production of some of the benchmarking results. 
This work was funded under the embedded CSE programme of the ARCHER UK National Supercomputing Service (http://www.archer.ac.uk), grant GPU-eCSE01-48, \emph{Achieving high-fidelity continuum modelling on GPUs with Nektar++}. DM acknowledges support from the Royal Academy of Engineering under their Research Chair scheme.

%\begin{itemize}
%    \item Diego for help with testing. 
%\end{itemize}

%% file: Appendix.tex
\appendix
\section{Architecture-aware implementation in Nektar++} 
\label{app:algo}

\subsection{MemoryRegion class}
A new \lstinline[style=C++Style]{MemoryRegion} 
class  has been introduced containing a dual pointer which is designed to ensure proper memory synchronisation for a contiguous block of memory with minimal data movement between host/device. Inside the class, a pair of Boolean variables are used to indicate where the latest data is located. On the user side, they are controlled by three memory access qualifiers, \lstinline[style=C++Style]!ReadOnly!, \lstinline[style=C++Style]!WriteOnly! and \lstinline[style=C++Style]!ReadWrite!. We provide further implementation details here to supplement the definitions given in Section~\ref{sec:implementation:mem}.
\begin{itemize}
	\item \textbf{ReadOnly}: The API returns an error message if the data has not been previously initialised by using a \lstinline[style=C++Style]!WriteOnly! memory access. Memory transfer occurs if the data in the current memory space was marked as invalid. A constant pointer is returned by the API and the data in the requested memory space is now marked as valid. The data in the other memory space remain valid.
	\item \textbf{WriteOnly}: Memory is allocated and marked as initialised, if not yet filled. No memory transfer occurs and a non-constant pointer is returned by the API. The data in the requested memory space is now marked as valid, while the data in the other memory space is marked as invalid.
	\item \textbf{ReadWrite}: The API returns an error message if the data has not been previously initialised by using a \lstinline[style=C++Style]!WriteOnly! memory access. Memory transfer occurs if the data on the requested memory space was marked as invalid and a non-constant pointer is returned by the API. The data in the requested memory space is now marked as valid, while the data in the other memory space is marked as invalid.
\end{itemize}

From the users' perspective, the memory can be accessed by simply calling \lstinline[style=C++Style]!GetPtr<MemSpace, MemAccess>!, where \lstinline[style=C++Style]!MemSpace! and \lstinline[style=C++Style]!MemAccess! are, respectively, the desired memory space and memory access qualifier, regardless where the live data is actually located. The memory within an instance of a \lstinline[style=C++Style]!MemoryRegion! class is allocated lazily the first time when a \lstinline[style=C++Style]!WriteOnly! access is called, and deallocated automatically when the variable goes out of its scope.

\subsection{Data structure and storage layout}
\begin{figure}[htb]
	\centering
	\includegraphics[width=0.6\textwidth]{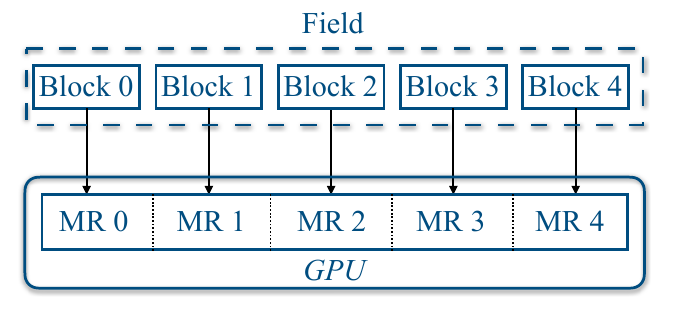}
	\caption{Field and Block data structure association. Field contains a series of blocks representing similar properties of elements that are associated with a contiguous block of memory on the device (i.e. GPU).}
	\label{fig:block_structure}
\end{figure}

Based on the properties of spectral element data, a Field/Block data structure is introduced to allow flexibility of implementation and kernel optimisation for performance-critical spectral element operators. We define a high-level class,  \lstinline[style=C++Style]!Field!, to represent a dataset associated with spectral elements within a parallel partition. The \lstinline[style=C++Style]!Field! class supports various data types through a template parameter \lstinline[style=C++Style]!TData! allowing for support for both single and double precision calculation. In Nektar++, there also exist two frequently-used concepts: \textit{coefficient} state, which represents the space of all expansion basis coefficients (or degrees of freedom, DOFs), and; \textit{physical} state, which represents the space of all quadrature points. Appropriate use of the \textit{coefficient} and \textit{physical} states is reinforced through an additional template parameter when declaring the \lstinline[style=C++Style]!Field! class.  The syntax is shown as follow;
\begin{lstlisting}[style=C++Style]
	// Initialize a Field object in coefficient space.
	auto FieldIn = Field<TData, FieldState::Coeff>(...);
	// Initialize a Field object in physical space.
	auto FieldOut = Field<TData, FieldState::Phys>(...);
\end{lstlisting}
The \lstinline[style=C++Style]!Field! class contains a vector of \lstinline[style=C++Style]!Block! objects. Each \lstinline[style=C++Style]!Block! 
includes an instance of a \lstinline[style=C++Style]!MemoryRegion! object and represents a collection of elements of the same "type" (same shape, basis, quadrature, and deformation type). A schematic representation of the Field/Block data structure is presented in Figure \ref{fig:block_structure}. 
Each individual \lstinline[style=C++Style]!Block! instance of a \lstinline[style=C++Style]!Field! object can then be accessed using an index-based loop of the form: 
\begin{lstlisting}[style=C++Style]
	for (unsigned int blk = 0; blk < field.GetBlocks().size(); ++blk)
	{
		// Fetch pointer from current block.
		auto ptr = field.GetBlocks()[blk].template GetPtr<DeviceSpace, WriteOnly>;
		...
		// Initialize data on a Device backend.
		...
	}
\end{lstlisting}

Conceptually, we can reshape each block into a multi-dimensional array. The default storage layout (from leading dimension to the last dimension) is quadrature points/modal coefficients, followed by elements, and finally components (i.e. components of a vector field), for example:
\begin{lstlisting}[style=C++Style]
	auto inptr  = inblock.template GetPtr<HostSpace, ReadOnly>;
	auto outptr = outblock.template GetPtr<HostSpace, WriteOnly>;
	for (unsigned int n = 0; n < block.GetNumComponents(); ++n)
	{
		for (size_t e = 0; e < block.GetNumElementsWithPadding(); ++e)
		{
			for (unsigned int p = 0; p < block.GetNumData(); ++p)
			{
				// Access each data point sequentially.
				auto point = *inptr;
				...
				// Write each data point sequentially.
				*outptr = ...
				
				// Increment pointers.
				inptr++;
				outptr++;
			}
		}
	}
	// Set block interleave parameter.
	outblock.template SetInterleaveWidth<TData>(1);
\end{lstlisting}
where \lstinline[style=C++Style]!block.GetNumData()! returns the number of degrees of freedom or quadrature points within one element. As previously mentioned, different algorithms may require different storage layouts for optimal performance. For example, for the AVX execution space, we cluster points from four to eight elements into the same SIMD register. Loading them takes minimum overhead if they are stored contiguously in a so-called \textit{interleaved} data layout. A sample code to access interleaved data is given below:
\begin{lstlisting}[style=C++Style]
	auto inptr  = inblock.template GetPtr<HostSpace, ReadOnly>;
	auto outptr = outblock.template GetPtr<HostSpace, WriteOnly>;
	for (unsigned int n = 0; n < block.GetNumComponents(); ++n)
	{
		for (size_t e = 0; e < block.GetNumElmtGroups(simd_type::width); ++e)
		{
			for (unsigned int p = 0; p < block.GetNumData(); ++p)
			{
				// Load data into simd_type.
				// inptr[0] : the p-th point of the (e * width + 0)-th element
				// inptr[1] : the p-th point of the (e * width + 1)-th element
				// inptr[2] : the p-th point of the (e * width + 2)-th element
				// ...
				simd_type point(inptr);
				...
				// Store data using simd_type.
				*((simd_type *)outptr) = ...
				
				// Increment pointers by vector width.
				inptr += simd_type::width;
				outptr += simd_type::width;
			}
		}
	}
	// Set block interleave parameter.
	outblock.template SetInterleaveWidth<TData>(simd_type::width);
\end{lstlisting}
where \lstinline[style=C++Style]!block.GetNumElmtGroups(simd_type::width)! returns the number of element (vector) groups associated with a SIMD vector width of \lstinline[style=C++Style]!simd_type::width! elements.

\subsubsection{Operator class hierarchy}
\begin{figure}
	\tiny
	\centering
	\begin{forest}
		for tree={
			font=\sffamily\bfseries,
			line width=1pt,
			draw=linecol,
			rrect,
			align=center,
			child anchor=north,
			parent anchor=south,
			l sep+=25pt,
			edge path={
				\noexpand\path[color=linecol, rounded corners=0pt, >={Latex[length=8pt]}, line width=1pt, <-, \forestoption{edge}]
				(!u.parent anchor) -- +(0,-20pt) -|
				(.child anchor)\forestoption{edge label};
			},
		}
		[Operator$\langle$TData$\rangle$, drop shadow, inner color=col1in, outer color=col1out, calign=child, calign child=1
		[ElmtOp$\langle$TFieldIn{,} TFieldOut{,} TData$\rangle$, drop shadow, inner color=col2in, outer color=col2out, calign=child, calign child=3
		[BwdTransOp, drop shadow, inner color=col3in, outer color=col3out] (a)
		[IProductWRTBaseOp, drop shadow, inner color=col3in, outer color=col3out]
		%[IProductWRTDerivBaseOp, drop shadow, inner color=col3in, outer color=col3out]
		[PhysDerivOp, drop shadow, inner color=col3in, outer color=col3out]
		[HelmholtzOp, drop shadow, inner color=col3in, outer color=col3out]
		]
		]
		\node[text width=20mm, text centered, draw=linecol, drop shadow, line width=1pt, font=\sffamily\bfseries,
		, rounded corners=4pt, inner color=col1in, outer color=col1out] (a0) at (-40mm, -0mm) {BlockOperator$\langle$TData$\rangle$};
		\node[text width=37.5mm, text centered, draw=linecol, drop shadow, line width=1pt, font=\sffamily\bfseries,
		, rounded corners=4pt, inner color=col2in, outer color=col2out] (a00) at (-40mm, -15mm) {ElmtBlockOp$\langle$TFieldIn{,} TFieldOut{,} TData$\rangle$};
		%\draw[color=linecol, line width=1.0pt, >={Latex[length=8pt]}, ->] (a00.north) -- ++ (0,0) |- (a0.south) ++ -- (0,-5mm) ;
		\draw[color=linecol, line width=1.0pt, >={Latex[length=8pt]}, ->] (a00.north) -| (a0.south);
		\node[text height=25mm, text width=60mm, draw=linecol, dashed] at (-22.5mm, -50mm) {};
		\node[text width=20mm, text centered, draw=linecol, inner color=col5in, outer color=col5out] (a1) at (-40mm, -40mm) {BwdTransBlockOp};
		\node[text width=20mm, text centered, draw=linecol, inner color=col5in, outer color=col5out] (b1) at (-40mm, -45mm) {BwdTransBlockOp};
		\node[text width=20mm, text centered, draw=linecol, inner color=col5in, outer color=col5out] (c1) at (-40mm, -50mm) {BwdTransBlockOp};
		\node[text width=20mm, text centered] (d1) at (-70mm, -55mm) {...};
		\node[text width=20mm, text centered, draw=linecol, inner color=col5in, outer color=col5out] (e1) at (-40mm, -60mm) {BwdTransBlockOp};
		\node[text width=27.5mm, text centered, draw=linecol, inner color=col5in, outer color=col5out] (a2) at (-10mm, -40mm) {BwdTransDeviceSumFac};
		\node[text width=27.5mm, text centered, draw=linecol, inner color=col5in, outer color=col5out] (b2) at (-10mm, -45mm) {BwdTransDeviceStdMat};
		\node[text width=27.5mm, text centered, draw=linecol, inner color=col5in, outer color=col5out] (c2) at (-10mm, -50mm) {BwdTransSerialStdMat};
		\node[text width=27.5mm, text centered] (d2) at (-40mm, -55mm) {...};
		\node[text width=27.5mm, text centered, draw=linecol, inner color=col5in, outer color=col5out] (e2) at (-10mm, -60mm) {BwdTransSerialSumFac};
		\draw[color=linecol, line width=1.0pt, >={Latex[length=8pt]}, ->] (a1.west) -- ++ (-12.5mm,0) |- (a00.west);
		\draw[color=linecol, line width=1.0pt, >={Latex[length=8pt]}, ->] (b1.west) -- ++ (-12.5mm,0) |- (a00.west);
		\draw[color=linecol, line width=1.0pt, >={Latex[length=8pt]}, ->] (c1.west) -- ++ (-12.5mm,0) |- (a00.west);
		\draw[color=linecol, line width=1.0pt, >={Latex[length=8pt]}, ->] (e1.west) -- ++ (-12.5mm,0) |- (a00.west);
		\draw[color=linecol, line width=1.0pt, >={Latex[length=5pt]}, ->] (a2) -- (a1);
		\draw[color=linecol, line width=1.0pt, >={Latex[length=5pt]}, ->] (b2) -- (b1);
		\draw[color=linecol, line width=1.0pt, >={Latex[length=5pt]}, ->] (c2) -- (c1);
		\draw[color=linecol, line width=1.0pt, >={Latex[length=5pt]}, ->] (e2) -- (e1);
		\draw[color=linecol, line width=1.0pt, >={Diamond[length=10pt]}, ->] (-35mm, -36.5mm) -- (-35mm, -31.5mm);
		\node[draw=none, inner color=white, outer color=white]  at (0mm, 5.5mm) {Operates on Field objects};
		\node[draw=none, inner color=white, outer color=white]  at (-40mm, 5mm) {Operates on Block objects};
		\node[draw=linecol, inner color=white, outer color=white] (e3) at (30mm, -60mm) {Specialised templated kernels};
		\draw[color=linecol, line width=0.5pt, >={Diamond[length=8pt]}, ->] (e3.west) -- ++ (-4mm,0) |- (a2.east);
		\draw[color=linecol, line width=0.5pt, >={Diamond[length=8pt]}, ->] (e3.west) -- ++ (-4mm,0) |- (b2.east);
		\draw[color=linecol, line width=0.5pt, >={Diamond[length=8pt]}, ->] (e3.west) -- ++ (-4mm,0) |- (c2.east);
		\draw[color=linecol, line width=0.5pt, >={Diamond[length=8pt]}, ->] (e3) -- (e2);
	\end{forest}
	\caption{Schematic representation of \lstinline[style=C++Style]!Operator! class hierarchy}
	\label{fig:operator_class}
\end{figure}
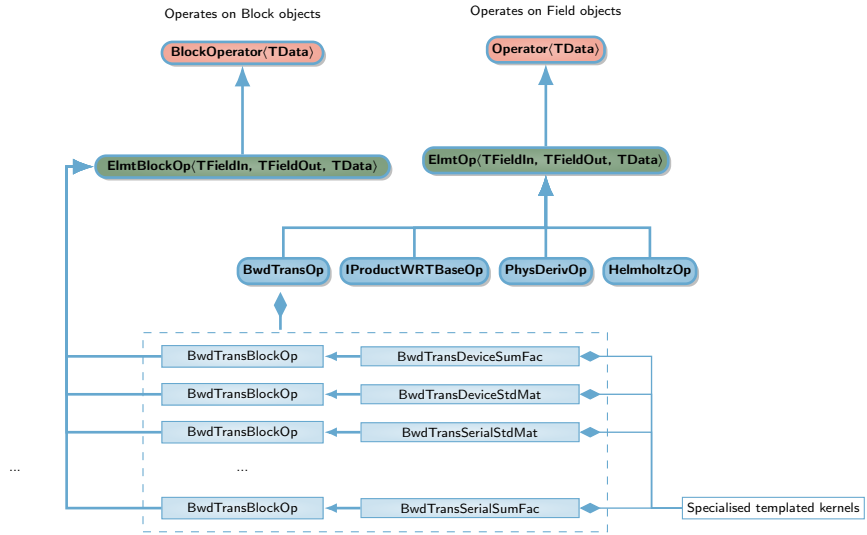

The Field/Block data structure described above allows the use of different implementation (i.e. StdMat, SumFac, or SumFacTOP) for each Block. For this purpose, a templated \lstinline[style=C++Style]!Operator! base class has been introduced, where each instance of an \lstinline[style=C++Style]!Operator! class contains a vector of instances of the \lstinline[style=C++Style]!BlockOperator! base class, allowing implementation specialisation for each Block. An example of implementation for the BwdTrans operator would take the following form: 
\begin{lstlisting}[style=C++Style]
	template <typename TData>
	class BwdTransOp : public ElmtOp<FieldState::Coeff, FieldState::Phys, TData>
	{
		void Apply(Field<TData, FieldState::Coeff> &in,
		Field<TData, FieldState::Phys> &out) override
		{
			// Loop over the blocks.
			for (unsigned int blk = 0; blk < this->m_blockOp.size(); ++blk)
			{
				// Fetch Block from Field objects.
				auto &inblock  = in.GetBlocks()[blk];
				auto &outblock = out.GetBlocks()[blk];
				
				// Use Block-specific implementation.
				this->m_blockOp[blk]->Apply(inblock, outblock);
			}
		}
		
		...
		
		// Vector of Block-specific BwdTransBlockOp implementations.
		std::vector<std::shared_ptr<BwdTransBlockOp<TData>>> m_blockOp;
	};
\end{lstlisting}
Specific operators, with different implementation or execution, can then be derived from the base class. Figure \ref{fig:operator_class} outlines the \lstinline[style=C++Style]!Operator! class hierarchy with a verbose expansion on BwdTrans operator.

\section{Additional results}
\subsection{Affine versus deformed elements on GH200}
\label{app:results:affine}

\begin{figure}[htb]
\includegraphics[width=\textwidth]{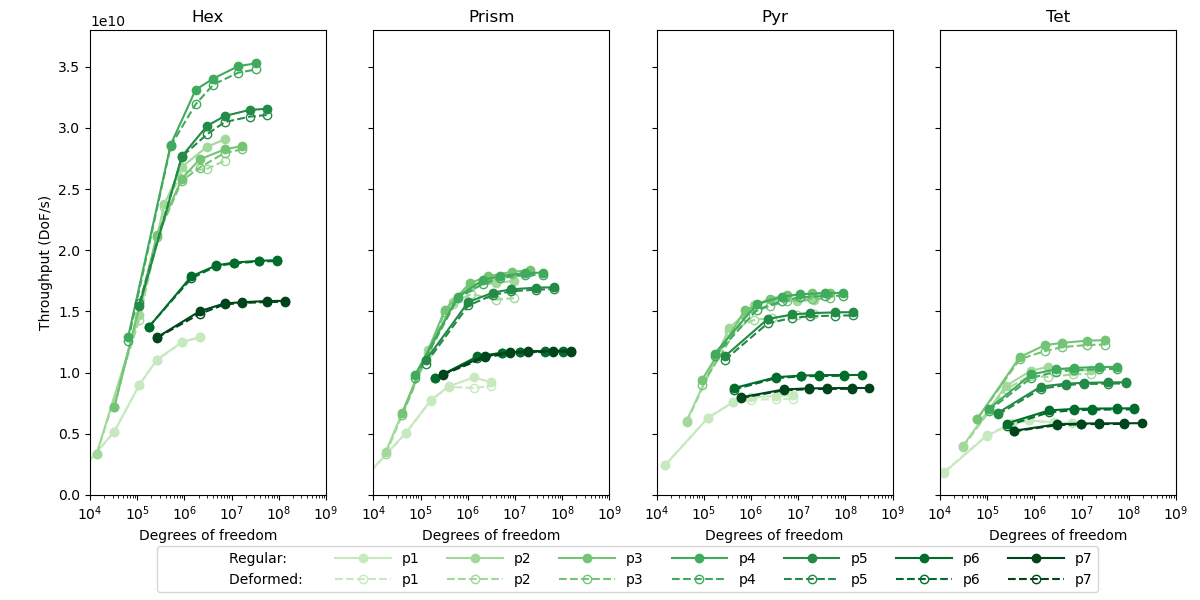}
\caption{NVIDIA GH200 Grace Hopper Superchip throughput performance versus elemental degrees of freedom on mass operator for hexahedral (Hex), prismatic (Prism), pyramidic (Pyr) and tetrahedral (Tet) elements using the sum-factorisation threaded on output-point (SumFacTOP) implementation. The solid lines indicate throughput assuming regular (Regular) elements that have an affine transformation where as the dashed lines indicate the throughput for curvilinear (Deformed) elements. }
\label{fig:reg_def_gh200_mass}
\end{figure}

Figure~\ref{fig:reg_def_gh200_mass} shows performance of mass operator on affine elements and curvilinear elements. Minimal differences are observed between the two cases.

\begin{figure}[htb]
\includegraphics[width=\textwidth]{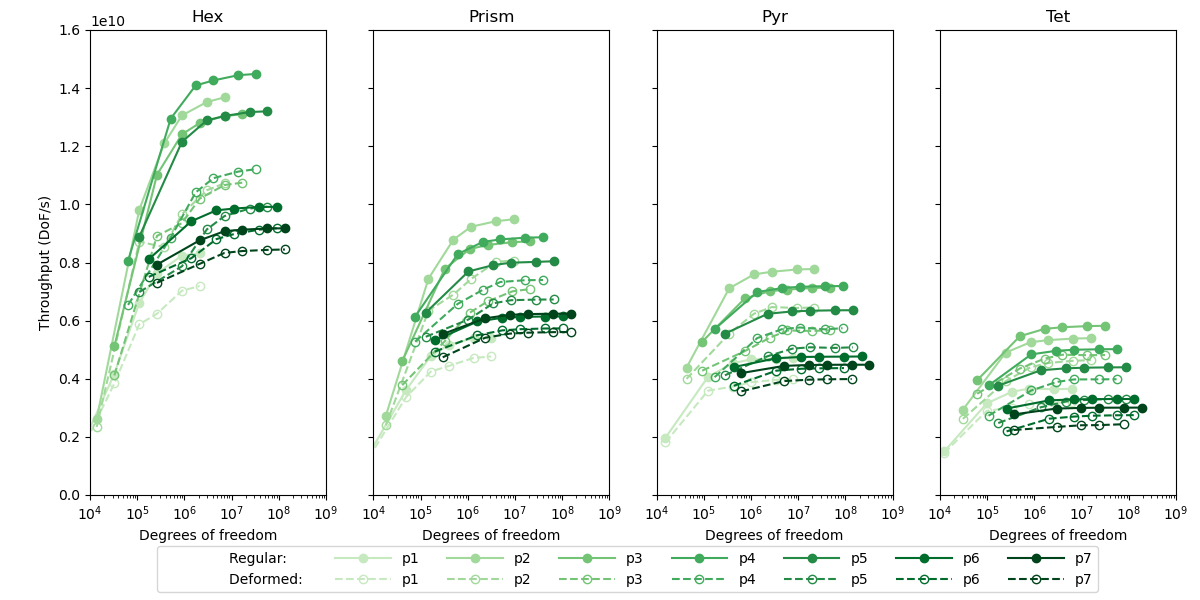}
\caption{NVIDIA GH200 Grace Hopper Superchip throughput performance versus elemental degrees of freedom on Helmholtz operator for hexahedral (Hex), prismatic (Prism), pyramidic (Pyr) and tetrahedral (Tet) elements using the sum-factorisation threaded on output-point (SumFacTOP) implementation. The solid lines indicate throughput assuming regular (Regular) elements that have an affine transformation where as the dashed lines indicate the throughput for curvilinear (Deformed) elements. }
\label{fig:reg_def_gh200_helm}
\end{figure}

Figure~\ref{fig:reg_def_gh200_helm} shows performance of Helmholtz operator on affine elements and curvilinear elements. This operator requires geometric factors at each quadrature point increasing data transfer costs. At low polynomial orders there is increased throughput when using regular affine elements; at higher polynomial orders, the benefit is limited.

\subsection{Performance on AMD EPYC 9554}
\label{app:results:epyc}

\begin{figure}[htb]
    \begin{center}
        \includegraphics[width=0.9\textwidth]{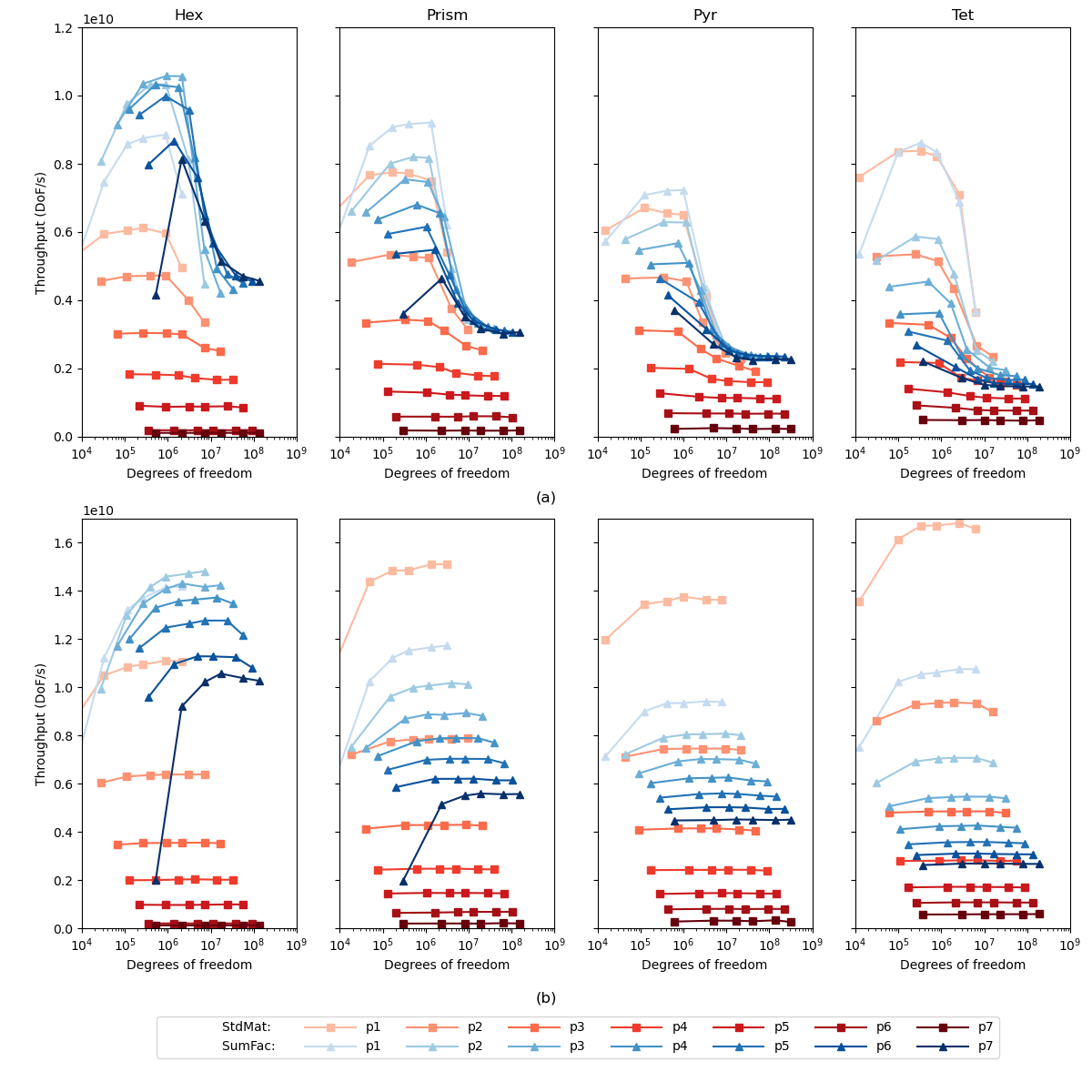}
    \end{center}
    \caption{Throughput results for the AMD EPYC 9554 processor for (a) deformed and (b) regular elements for the Helmholtz operator on all element types using the SumFac and StdMat implementation strategies.}
    \label{fig:epyc-9554}
\end{figure}

Figure~\ref{fig:epyc-9554} highlights throughput results on a dual-socket AMD EPYC 9554 processor, which has a significantly increased number of cores (64 per socket) and L3 cache per core (256MB total, 4MB/core) compared to the Intel Xeon 6526Y (16 per socket, 37.5MB total or 2.3MB/core).

